\documentclass[a4paper,12pt]{article}

\usepackage[a4paper, centering]{geometry}
\usepackage{longtable}
\usepackage{fancyhdr}

\usepackage{cancel}
\usepackage{enumitem}

\setitemize{noitemsep,topsep=0pt,parsep=0pt,partopsep=0pt,leftmargin=*}

\usepackage{amsmath}
\usepackage{amsfonts}
\usepackage{amsthm}
\usepackage{amssymb}
\usepackage{upref}
\usepackage{color}
\usepackage[table,xcdraw,svgnames]{xcolor}
\usepackage{graphicx}

\usepackage[colorlinks]{hyperref}
\hypersetup{
	allbordercolors=.,
	citecolor=NavyBlue,
	linkcolor=DarkRed,
	urlcolor=NavyBlue
}

\usepackage{array}
\usepackage[utf8]{inputenc}
\usepackage{algorithm,algorithmicx,algpseudocode}

\graphicspath{{./pics_all/}}

\theoremstyle{plain}
\newtheorem{thm}{Theorem}[section]

\newtheorem{prop}[thm]{Proposition}

\newtheorem{rem}[thm]{Remark}
\newtheorem{lemma}[thm]{Lemma}
\theoremstyle{definition}
\newcommand{\Per}{\operatorname{Per}}

\newcommand{\di}{\operatorname{div}}
\newcommand{\bo}[1]{{\bf#1}}

\usepackage{indentfirst} 
\newtheorem{problem}{Problem}

\title{Numerical shape optimization among convex sets}

\author{Beniamin Bogosel\footnote{CMAP, Ecole Polytechnique, route de Saclay, 91128 Palaiseau, France (\nolinkurl{beniamin.bogosel@polytechnique.edu})}}

\begin{document}
	
	\maketitle
	
	\begin{abstract}
		This article proposes a new discrete framework for approximating solutions to shape optimization problems under convexity constraints. The numerical method, based on the support function or the gauge function, is guaranteed to generate discrete convex shapes and is easily implementable using standard optimization software. The framework can handle various objective functions ranging from geometric quantities to functionals depending on partial differential equations. Width or diameter constraints are handled using the support function. Functionals depending on a convex body and its polar body can be handled using a unified framework. 
	\end{abstract}



\section{Introduction}
Shape optimization problems involve the minimization or maximization of functionals having geometric shapes as variables. Given a shape optimization problem, multiple non-trivial aspects can be investigated: existence of optimal shapes, study of optimality conditions, regularity and qualitative properties of the optimal shapes. The monographs \cite{henrot-pierre-english}, \cite{delfour-zolesio} give an overview of challenging aspects and methods used in this field. Cases where the optimal shape can be explicitly identified are quite rare. When dealing with convex sets the proof of existence of optimal shapes is often straightforward, due to the classical Blaschke selection theorem \cite[Theorem 1.8.7]{schneider}. However, as the recent works \cite{LNconv09,jimmy_2021} show, studying optimality conditions under convexity constraints may be quite challenging. The optimal shapes for the problems considered are not explicitly known, in general. This motivates the development of numerical methods for approximating solutions to shape optimization problems among convex shapes. 

The convexity constraint poses difficulties in the numerical implementation since it restricts the class of admissible domain perturbations. The recent work \cite{bartels_wachsmuth} shows how to handle the convexity constraint by working with deformed meshes and constraining the admissible deformations. In \cite{BHcw12} and \cite{ABcw} the authors use a truncated spectral decomposition of the support function and handle a wide variety of constraints. One drawback of using truncated spectral decompositions for the support function is the smoothness of the support function in the discrete setting. It is well known \cite[Section 1.7]{schneider} that smooth support function correspond to strictly convex shapes. In particular, segments in the boundary are only captured in an approximate way using truncations of spectral decompositions. This is a fundamental aspect, since imposing a convexity constraint can naturally produce optimal shapes having segments in the boundary. In the paper \cite{steklov_diam} the authors use a finite difference method for parametrizing the support function, allowing discontinuities in the first derivative and capturing efficiently segments in the boundary of the convex domains. 
In this article we further develop ideas that help solve numerically shape optimization problems for convex sets, including discontinuities in the derivative of the support function, corresponding to segments in the boundary of the optimal shape. For the discretization of the support funciton a more rigorous method, compared to \cite{steklov_diam}, is proposed and studied in detail. In addition, a numerical framework using the gauge function is also proposed. Various applications are presented and the code used for producing the numerical results is freely available. 

Other works in the literature deal with numerical aspects related to the convexity constraint. We mention \cite{oudet-convex-hull} where the authors propose a parametrization using supporting half-spaces.  In \cite{merigot-oudet} the discretization of optimization problems with convexity constraints is investigated. In \cite{galerkin-convex} the author proposes a Galerkin approximation theory for convex sets.

\begin{figure}
	\centering
	\includegraphics[width=0.4\textwidth]{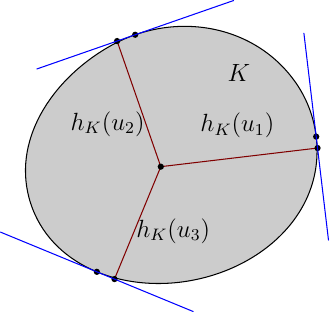}
	\caption{Illustration of the support function of a convex body.}
	\label{fig:suppfunc}
\end{figure}

A convex body $K \subset \Bbb{R}^d$ is a compact convex set with non-void interior. The support function $h_K: \Bbb{S}^{d-1} \to \Bbb{R}$ of a convex body $K\subset \Bbb{R}^d$ is defined by
\begin{equation}
 h_K(u) = \max_{x \in K} (x \cdot u) 
\label{eq:def_supp_func}
\end{equation}
or alternatively, $h_K(u)$ is the distance from the origin to the supporting plane orthogonal to the direction $u \in \Bbb{S}^{d-1}$. An illustration for the two dimensional case is given in Figure \ref{fig:suppfunc}. It is also possible to define the support function on the whole space $H_K:\Bbb{R}^d \to \Bbb{R}$ by 
\begin{equation}
 H_K(u) = \max_{x \in K} x \cdot u,
 \label{eq:def_supp_func_extended}
\end{equation}
and note that $h_K = H_K(u/|u|)$ when $u\neq 0$. In other words, $H_K$ is the positive $1$-homogeneous function which coincides with $h_K$ on the unit sphere. The concept of support function is classical in convex geometry and the reader can consult \cite{schneider} for more details and properties. In Section \ref{sec:suppfunc} all aspects of support functions that are relevant to this work are recalled.

The definition above shows, in particular, that the support function is well adapted for dealing numerically with width or diameter constraints. This was already observed in the previous works \cite{ABcw}, \cite{BHcw12}, \cite{BHL17} or \cite{steklov_diam}. In \cite[Section 1.7]{schneider} it is shown that the support function $H_K$ in \eqref{eq:def_supp_func_extended} is sublinear (positive $1$-homogeneous verifying $f(u+v)\leq f(u)+f(v),\ \forall u,v \in \Bbb{R}^d$) and for every sublinear function $f:\Bbb{R}^d \to \Bbb{R}$  there exists a convex body $K \subset \Bbb{R}^d$ such that $H_K = f$. In dimension two, in the case $K$ is of class $C^1$, a necessary and sufficient condition for $h_K \in C^1$ to be a support function of a convex body is to verify $h_K(\theta)+h_K''(\theta) \geq 0$ (in the sense of distributions) for all $\theta \in [0,2\pi]$ (see \cite[Chapter 1]{schneider} or \cite{BHcw12} for example). In higher dimensions the characterization of the constraint becomes more complex as shown in \cite{ABcw}, for the three dimensional case. 

\begin{figure}
	\centering
	\begin{tabular}{ccc}
	\includegraphics[height=0.2\textwidth]{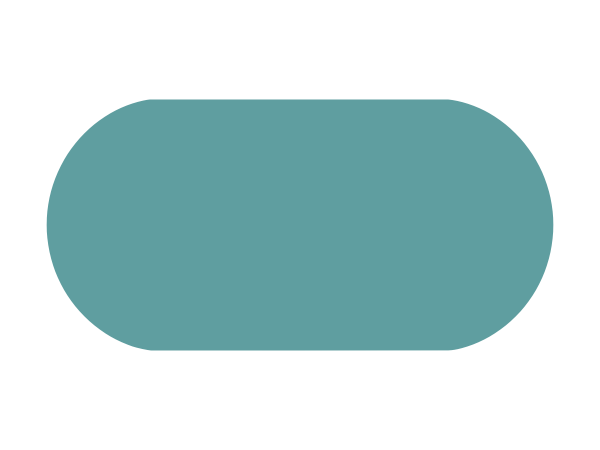}&
	\includegraphics[height=0.2\textwidth]{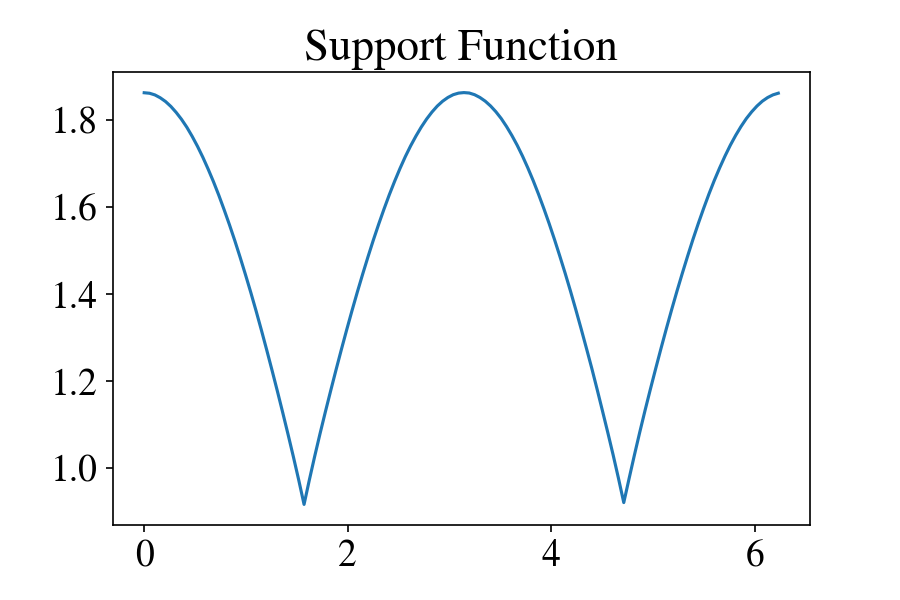}&
	\includegraphics[height=0.2\textwidth]{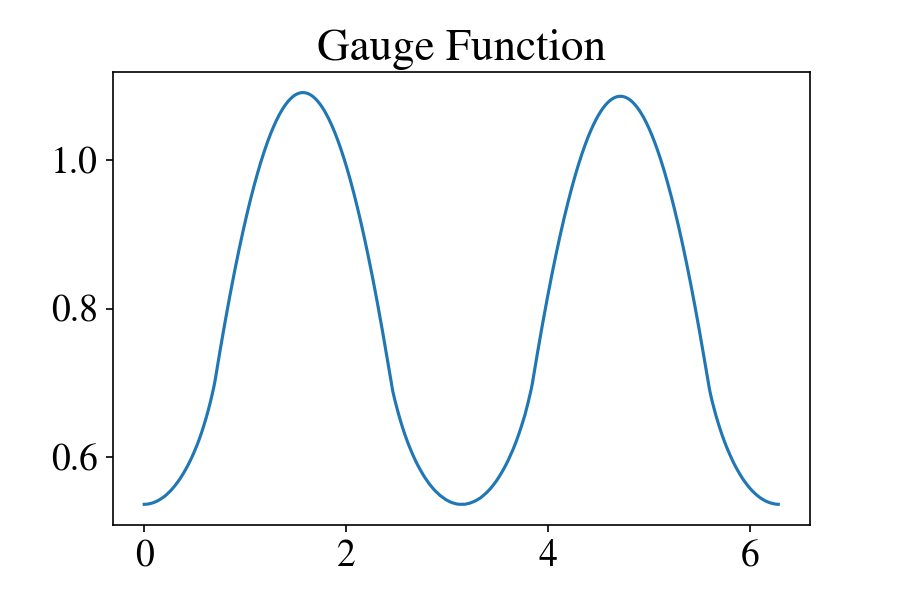}\\
	\includegraphics[height=0.2\textwidth]{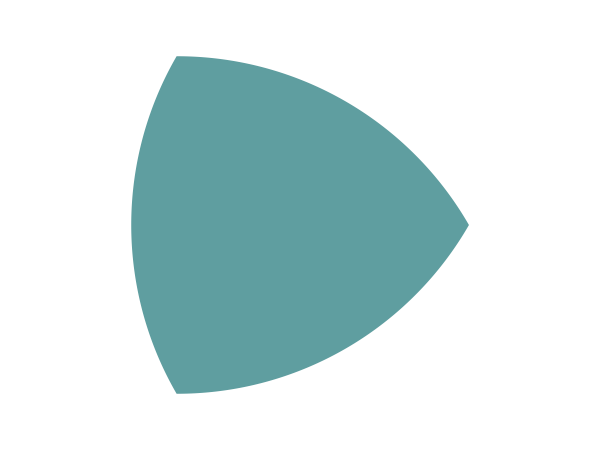}&
	\includegraphics[height=0.2\textwidth]{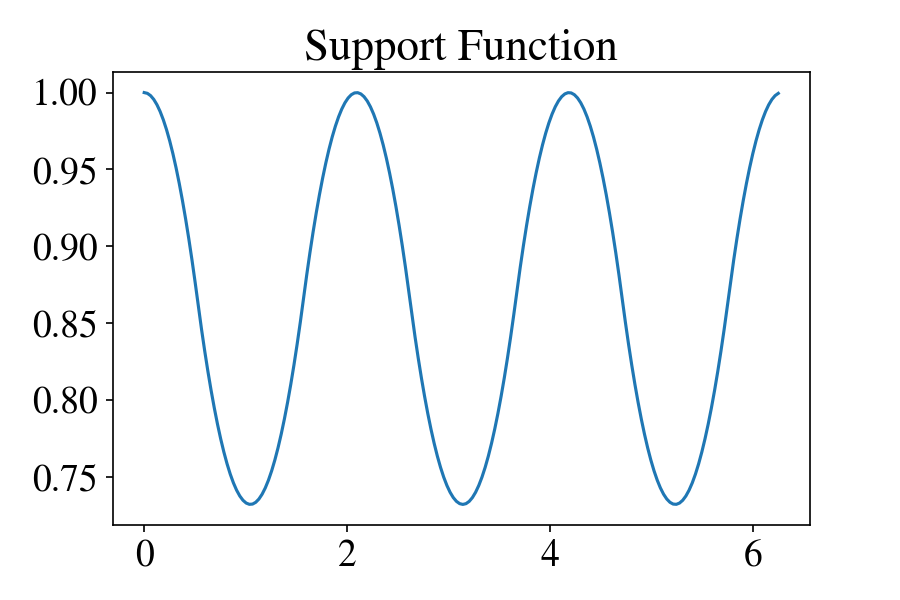}&
	\includegraphics[height=0.2\textwidth]{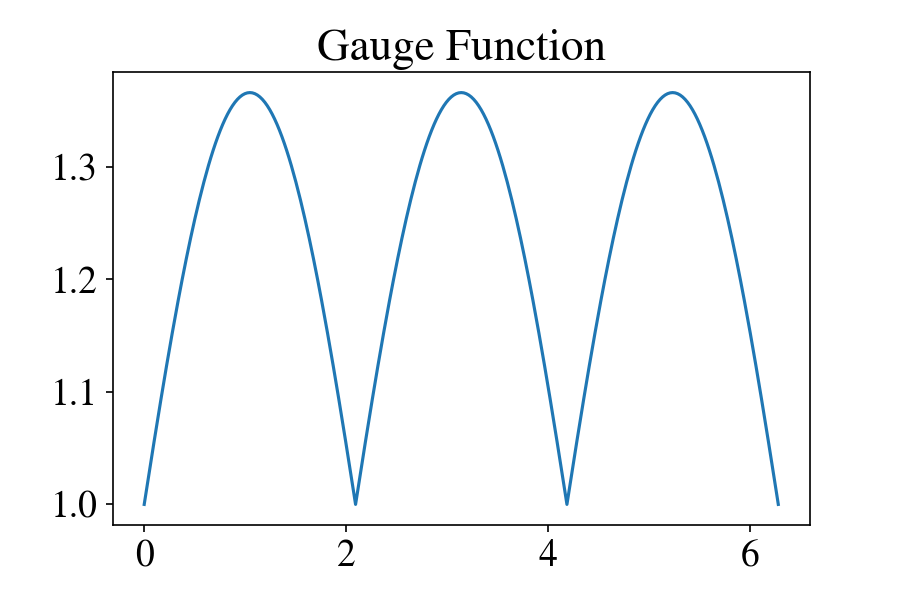}
\end{tabular} 
	\caption{Two examples of convex shapes together with their support and gauge functions. Segments in the boundary correspond to singular points for the support function. Angular points generate singularities for the gauge function.}
	\label{fig:supp_examples}
\end{figure}

In Figure \ref{fig:supp_examples} two examples of shapes and their associated support functions are shown, for the classical Reuleaux triangle and a stadium like shape. It can clearly be seen that segments in the boundary of the stadium correspond to points where the support function has a discontinuous derivative. In \cite[Cor 1.7.3]{schneider} it is shown that at all points where the supporting plane intersects $K$ at exactly one point, the support function is differentiable. In particular, segments in the boundary of a two dimensional convex domain produce discontinuities in the first derivative of the associated support function. Therefore, the discretization of the support function should allow such discontinuities in the derivative in order to capture segments in the boundary.
In \cite{steklov_diam} a numerical method was introduced which does not impose any regularity assumptions on the discretization of the support function. The method consists in considering the values of the support function in dimension two at a uniform discretization of the unit circle. All differential quantities involved in the computations are approximated using finite differences. In this paper we propose an alternative method, which in addition to being consistent, is guaranteed to produce discrete convex shapes in the optimization process. From a numerical point of view, the proposed numerical method has the same complexity as the one shown in \cite{steklov_diam}.

Another natural parametrization of convex sets can be achieved using radial functions with respect to an interior point. In dimension two it turns out that the inverse of the radial distance to a fixed origin has properties which allows again the use of efficient numerical methods for discretizing convex shapes. Given $K\subset \Bbb{R}^d$ a convex set containing the origin, consider $\rho_K : [0,2\pi] \to \Bbb{R}_+$ to be a radial function for $K$. In other words, $\theta \mapsto \rho_K(\theta)\begin{pmatrix}
\cos \theta\\
\sin \theta
\end{pmatrix}$ is a parametrization for $\partial K$. This allows us to define the associated gauge function $\gamma_K :[0,2\pi]\to \Bbb{R}_+$ by
\begin{equation}
  \gamma_K(\theta) = 1/\rho_K(\theta).
  \label{eq:gauge}
\end{equation}
The gauge function is related to the support function via the polar body. Given a convex body $K$, the polar body is $K^\circ = \{y \in \Bbb{R}^d : x\cdot y \leq 1, \forall x \in K\}$. The gauge function of the body $K$ is equal to the support function for $K^\circ$. See \cite{schneider} for more details. Compared with the support function, the gauge function has singularities at corners (instead of segments). Two examples are shown in Figure \ref{fig:supp_examples} illustrating the smoothness of the gauge function for segments in the boundary and singularities coming from angular points. Details regarding the usage of the gauge function in numerical simulations are presented in Section \ref{sec:gauge}.

The paper is structured as follows. In Section \ref{sec:suppfunc} we describe the discretization of the support function and we analyze the behavior of the resulting numerical method. In Section \ref{sec:gradient} we describe how the gradient of a general shape functional can be computed in terms of the parameters describing the support function. Section \ref{sec:applic-supp} illustrates how the proposed discretization using the support function approximates solutions to various shape optimization problems. Section \ref{sec:gauge} presents the discrete framework based on the gauge function and shows a few applications. 

The main purpose of this paper is to provide new, rigorous and flexible numerical tools for studying shape optimization problems among convex sets. In addition to the presentation of the numerical methods, various applications are shown. Many of these are classical, nevertheless, here are some new results that may motivate further theoretical study:
\begin{itemize}
	\item For $1\leq k\leq 10$, the shape which maximizes the $k$-th eigenvalue of the Dirichlet-Laplace operator among constant width shapes is the Reuleaux triangle. (Problem \ref{pb:max-lamk-const-width})
	\item For $1\leq k\leq 10$, the shape which maximizes the $k$-th eigenvalue of the Dirichlet-Laplace operator among convex shapes having a given minimal width is the equilateral triangle. (Problem \ref{pb:max-lamk-width})
\end{itemize}
Some of the codes that were used to produce the numerical results in the paper are provided:
\begin{center}
\href{https://github.com/bbogo/ConvexSets}{\nolinkurl{https://github.com/bbogo/ConvexSets}}
\end{center}
These codes can be used to approximate solutions to other shape optimization problems among convex sets by simply changing the evaluation of the objective function and the associated shape derivative, as indicated in Section \ref{sec:gradient}.

\section{Support function and its discretization}
\label{sec:suppfunc}
For a convex compact body $K\subset \Bbb{R}^2$ consider the associated support function $h_K$ defined by \eqref{eq:def_supp_func}. Let us briefly recall some of the basic properties of the support function. For a complete exposition with proofs the reader can consult \cite[Section 1.7]{schneider}. The support function $h_K$ can be identified with a continuous $2\pi$ periodic function $p: [0,2\pi] \to \Bbb{R}$. For $\theta\in [0,2\pi]$ denote the associated normal and tangential vectors $\bo r(\theta) = (\cos \theta,\sin\theta)$ and $\bo t(\theta) = (-\sin\theta,\cos\theta)$.  the set $H(K,\theta) = \{ x \in \Bbb{R}^2 : x\cdot \bo r(\theta) = p(\theta)\}$ is called the support line of $K$ at $\theta$. The set $F(K,\theta) = K\cap H(K,\theta)$ is called the support set of $K$ at $\theta$. In \cite[Section 1.7]{schneider} it is shown that $p$ is differentiable at $\theta$ if and only if the associated support set $F(K,\theta)$ contains only one point. An immediate consequence is that segments in the boundary of $K$ correspond to parameters $\theta$ where the support function has a discontinuity in the first derivative. In particular, constant width shapes have support functions at least of class $C^1$, since they do not contain non-trivial segments in their boundaries.

It is classical that given the support function $p$ of a strictly convex shape $K\subset \Bbb{R}^2$, $p$ is of class $C^1$ and a parametric representation of $\partial K$ is given by $\bo x(\theta) = p(\theta) \bo r(\theta)+p'(\theta) \bo  t(\theta)$ or, more explicitly,
\begin{equation} \begin{cases}
x_1(\theta) = p(\theta) \cos \theta - p'(\theta) \sin \theta \\
x_2(\theta) = p(\theta) \sin \theta + p'(\theta) \cos \theta.
\end{cases}
\label{eq:param}
\end{equation}
It is straightforward to see that at points where the support function $p$ is smooth, the tangent vector at the curve with the parametrization \eqref{eq:param} is given by $(p(\theta)+p''(\theta))\bo t(\theta)$. The parametrization \eqref{eq:param} is well defined at smooth points for $p$ provided $p''+p \geq 0$. In \cite{BHcw12}, the convexity of $K$ is characterized by 
\begin{equation}
p(\theta)+p''(\theta)\geq 0 \text{ for every } \theta \in [0,2\pi].
\label{eq:conv-constraint}
\end{equation}At points where $p$ is smooth, the condition above is clear. The constraint may also be interpreted in the sense of distributions, including points where $p''$ is not defined. Moreover, if $p$ satisfies $p+p''\geq 0$ in the sense of distributions then $p$ is the support function of a convex body $K$. The quantity $\varrho = p+p''$ is the curvature radius of $\partial K$. In particular, at regular points $\bo x(\theta)$ of $\partial K$ the curvature is given by $1/(p(\theta)+p''(\theta))$. In the following, we use the notation $\theta(\bo x)$ to represent the angle in $[0,2\pi]$ associated to the point $\bo x  \in \partial K$. Equivalently, $\theta(\bo x)$ is the angle the normal at $\bo x \in \partial K$ makes with the positive $x_1$-axis.

Suppose $\theta_0$ is an isolated point where $p'$ is not defined, corresponding to a segment $S_0$ in $\partial K$. Then the endpoints of the segment $S_0$ can be identified by
\[ \bo q_-=\lim_{\theta \to {\theta_0}, \theta<\theta_0} \bo x(\theta),\ 
\bo q_+=\lim_{\theta \to {\theta_0}, \theta>\theta_0} \bo x(\theta), \]
where $\bo x(\theta) = (x_1(\theta),x_2(\theta))$ is given by \eqref{eq:param}. 
It is straightforward then that the length of $S_0$ is given by $p'(\theta_+)-p'(\theta_-)>0$ (using the usual notation for half limits). This shows, in particular, that such a singularity for $p'$ cannot correspond to a local maximum of $p$.

The support function is particularly useful when dealing with numerical aspects related to the convexity constraint in shape optimization. In the works \cite{antunesMM16}, \cite{ABcw}, \cite{BHcw12} the authors use truncated Fourier series to parametrize support functions of convex sets. As already underlined in the previous paragraphs, this excludes the possibility of having segments in the boundary. Therefore, alternate discretization options are needed to handle such cases which arise quite often in practical situations. In view of the discussion above, discontinuities in the first derivative should be allowed by the choice of the discretization.

\subsection{Convexity constraint -- finite differences.} 

In view of the expression for parametrization \eqref{eq:param} and the continuous convexity constraint, any discretization strategy depends on the choice of the approximations of $p'$ and $p''$. For $N\geq 3$, consider an equidistant partition of $[0,2\pi]$ given by $\theta_j = 2\pi j/N$, $j=0,...,N-1$. The indices are considered periodic modulo $N$ whenever necessary. Denote $h=2\pi/N$ the distance between the points in the discretization and $p_j$ approximations of values $p(\theta_j)$ of the support function taken at $\theta_j$. The values $\bo p = (p_j)_{j=0}^{N-1}$ will be the \emph{optimization variables} for the problems presented in the first part of this work.

A first option, proposed in \cite{steklov_diam} is to approximate first derivatives using centered finite differences 
\begin{equation}
p(\theta_i) \approx (p_{i+1}-p_{i-1})/(2h),\ \  i=0,...,N-1.
\label{eq:finite-differences-order1}
\end{equation} 
At points where $p$ is of class $C^2$, finite differences can be used for approximating $p''$ giving the following inequalities characterizing the discrete convexity constraint:
\begin{equation}
 p_i+\frac{p_{i+1}-2p_i+p_{i-1}}{h^2} \geq 0,\ i=0,...,N-1.
 \label{eq:convexity-variant1}
\end{equation}

On the other hand, if $p$ is not smooth at $\theta_i$ then inequality \eqref{eq:conv-constraint} simply states that $\frac{p_{i+1}-p_i}{h}-\frac{p_i-p_{i-1}}{h} \geq -hp_i$. For $h$ small, supposing that $p_i$ is bounded from above, this inequality characterizes up to the first order in $h$ the positivity of the length of the segment associated to the singularity at $\theta_i$.

As a consequence, constraint \eqref{eq:conv-constraint} characterizes up to higher order terms in $h$ both the smooth and the non-smooth aspects of the convexity constraint from the general case. From a practical point of view \eqref{eq:conv-constraint} can be formulated as a set of $N$ linear inequality constraints, and easily be incorporated in standard optimization algorithms.

\begin{rem}
	If one chooses to parametrize the support function using piecewise affine functions on intervals $[\theta_i,\theta_{i+1}]$ then the constraint \eqref{eq:conv-constraint} implies that
	\[ \int_0^{2\pi} (p+p'')\varphi \geq 0 \text{ for any } \varphi \geq 0.\]
	Suppose $\psi_i$, $i=0,...,N-1$ is the classical basis for ${\bf P_1}$ finite elements. Then expressing the inequality above using the decompositions of $p$ and $\varphi$ in this basis one gets that the coefficients $\bo p = (p_i)$ of $p$, given again by $p_i=p(\theta_i)$ verify
	\begin{equation} \bo p M \bo \varphi - \bo p K \varphi \geq 0, \ \forall \varphi \geq 0.
	\label{eq:convex-finitel-2d}
	\end{equation}
	The matrices $M=(m_{ij})_{i,j=0}^{N-1}$ and $K= (k_{ij})_{i,j=0}^{N-1}$ are the mass and, respectively, rigidity matrices defined by 
	\[ m_{ij} = \int_0^{2\pi} \psi_i \psi_j,\ \ k_{ij} = \int_0^{2\pi} \psi'_i \psi_j', \ 0\leq i,j \leq N-1.\]
	In the case of ${\bf P_1}$ finite elements with periodic boundary conditions, using equidistant intervals $[\theta_i,\theta_{i+1}]$ leads to the inequalities
	\begin{equation}
	\frac{2}{3}p_i+\frac{1}{6}(p_{i+1}+p_{i-1})+\frac{p_{i+1}-2p_i+p_{i-1}}{h^2} \geq 0, \text{ for every }i=0,...,N-1.
	\label{eq:convexity-variant2}
	\end{equation}
	One can immediately note the similarity with \eqref{eq:convexity-variant1} and the fact that the two inequalities are equivalent up to first order in $h$.
	
	Nevertheless, any non-constant ${\bf P_1}$ support function admits local maxima where the derivative is discontinuous. As a consequence, no such function is the support function of a convex set. 
\end{rem}

\begin{figure}
	\centering
	\includegraphics[width=0.5\textwidth]{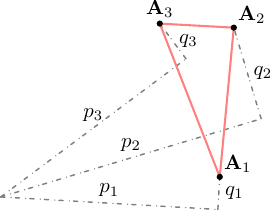}
	\caption{Construction of the triangle $\bo A_{1}\bo A_2 \bo A_{3}$ given by \eqref{eq:formula-tri}}.
	\label{fig:supp-tri}
\end{figure}

In the following we investigate the convexity properties of the discrete shape obtained using the proposed discretization of the support function. 
As stated before, we consider $p_i$ the values of the support function evaluated at $\theta_i$ and $q_i$ the chosen values for the approximation of $p'(\theta_i)$. Without loss of generality, we look at points for $\theta_1,\theta_2,\theta_3$ of the discretization. Then if 
\begin{equation}
\bo r_i = (\cos \theta_i,\sin \theta_i), \ \  \bo t_i = (-\sin \theta_i,\cos \theta_i)
\label{eq:radial-tangent}
\end{equation} are the normal and tangent vectors corresponding to $\theta_i$, the boundary points given by the parametrization \eqref{eq:param} of the support function are given by
\begin{equation}
\bo A_i = p_i \bo r_i +q_i \bo t_i, i=1,2,3.
\label{eq:formula-tri}
\end{equation}
An example is shown in Figure \ref{fig:supp-tri}.
The convexity property in the discrete setting amounts to the positivity of the oriented area of the triangle $\Delta \bo A_1\bo A_2\bo A_3$ for the points given in \eqref{eq:formula-tri}. The oriented area (in the trigonometric sense) in terms of the coordinates $(A_i^1,A_i^2)$ of $\bo A_i$, $i=1,2,3$, is obtained with the cross product formula
\begin{equation} \text{Area} (\Delta\bo  A_1\bo A_2\bo A_3) =\frac{1}{2}\left( (A_2^1-A_1^1)(A_3^2-A_2^2)-(A_3^1-A_2^1)(A_2^2-A_1^2)\right).
\label{eq:area-tri}
\end{equation}
Choosing $q_i = (p_{i+1}-p_{i-1})/(2h)$ corresponding to the approximation of the derivatives using centered finite differences \eqref{eq:finite-differences-order1} the area given in \eqref{eq:area-tri} can be expressed in terms of $p_i$. However, since we impose constraints on the discrete curvature radii defined by
\begin{equation}
\varrho_i = \frac{p_{i+1}-2p_i +p_{i-1}}{h^2}+p_i \geq 0
\label{eq:discrete-conv-constr}
\end{equation}
it is more convenient to express this area in terms of $\varrho_i, i=1,2,3$. In particular, we obtain
\begin{equation}\text{Area} (\Delta \bo  A_1\bo A_2\bo A_3) = \frac{1}{48}\left[6\varrho_1\varrho_2+12\varrho_1\varrho_3+6\varrho_2\varrho_3+(p_1-p_3)(\varrho_1-\varrho_3)\right] h^3+O(h^5),
\label{eq:expansion-area-fd}
\end{equation}
as $h\to 0$. The Mathematica script performing the symbolic computation is given in the Appendix. Observing that $\varrho_1-\varrho_3 = p_0-p_4+\left(1-2/h^2\right)(p_1-p_3)$, the coefficient of the leading term $h^3$ is given by 
\[ 6\varrho_1\varrho_2+12\varrho_1\varrho_3+6\varrho_2\varrho_3+(p_1-p_3)(p_0-p_4)+(1-2/h^2)(p_1-p_3)^2.\]
Therefore, in general, the positivity of $\varrho_1,\varrho_2,\varrho_3$ does not imply the positivity of the area of \eqref{eq:area-tri} as $h\to 0$. As a consequence, although this method was successfully used in \cite{steklov_diam}, the resulting discrete shapes may not be convex for $h$ small.

Another drawback that the classical finite differences discretization of the parametrization \eqref{eq:param} and of the convexity constraint \eqref{eq:convexity-variant1} are not invariant with respect to translations of the domain. This aspect is the key observation that allows us to propose an alternative discretization procedure which is consistent and which produces convex discrete shapes.

\subsection{Rigorous discrete convexity condition}

It is well known that if $p$ is the support function of the convex body $K\subset \Bbb{R}^2$ then $\overline p = p+a\cos \theta+b\sin \theta$ is the support function of the translated body $(a,b)+K$. This is a straightforward consequence of definition or of the parametrization \eqref{eq:param}. One may note that the classical centered finite differences are not exact when considering discretizations of translations given by $\overline p_i = p_i+a\cos \theta_i + b \sin \theta_i$. In order to remedy this we propose the following choices for approximating the first derivatives $p'(\theta_i)$ and the discrete curvature radii $p(\theta_i)+p''(\theta_i)$
\begin{eqnarray}
p'(\theta_i) \approx \frac{p_{i+1}-p_{i-1}}{2\sin h} , \ \ p(\theta_i)+p''(\theta_i) \approx \varrho_i= p_i+ \frac{p_{i+1}-2p_i+p_{i+1}}{2-2\cos h}.
\label{eq:new-discretization}
\end{eqnarray}
The observations below show that this discretization choice for first derivatives and curvature radii has multiple advantages:
\begin{itemize}
	\item As $h\to 0$ we have $\sin h = h+O(h^3)$, $2-2\cos h = h^2+O(h^4)$. Therefore, when $h \to 0$, at points where $p$ is smooth, the discretizations proposed in \eqref{eq:new-discretization} converge to $p'(\theta_i)$ and $p+p''(\theta_i)$, respectively. Therefore, the proposed discretization is {\bf consistent}.
	\item Formulas \eqref{eq:new-discretization} are linear and they are exact for support functions of the form $p(\theta) = c+a\sin \theta+b\cos \theta$. A first consequence is that the {\bf discretization process commutes with translations}: the numerical representations of two translated convex bodies are related by the same translation. Secondly, translated discretized convex bodies have the same discrete curvature radii given by \eqref{eq:new-discretization}.
	\item Using $q_i = (p_{i+1}-p_{i-1})/(2\sin h)$ in \eqref{eq:formula-tri} and computing the area of the triangle given in \eqref{eq:area-tri} we obtain
	\begin{equation}
	\text{Area}(\Delta \bo A_1\bo A_2\bo A_3) = [\varrho_2(\varrho_1+\varrho_3)+2\varrho_1\varrho_3\cos h]\sin^2(h/2)\tan(h/2).
	\label{eq:area-tri-rigorous}
	\end{equation}
	The Mathematica script performing the symbolic computation is given in the Appendix.
	Assuming the discrete radii of curvature are non-negative implies that the area of $\Delta \bo A_1\bo A_2\bo A_3$ is non-negative. Therefore, {\bf discrete shapes} constructed using \eqref{eq:new-discretization} in the parametrization \eqref{eq:param} {\bf are convex}, provided
  \begin{equation}
   \varrho_i=p_i+ \frac{p_{i+1}-2p_i+p_{i+1}}{2-2\cos h}=\frac{p_{i+1}+p_{i-1}-2p_i\cos h }{2-2\cos h} \geq 0, \ \ i=0,...,N-1.
   \label{eq:new-conv-constraint}
  \end{equation}
\end{itemize}
The constraints \eqref{eq:new-conv-constraint} are linear in the variables $(p_i)_{i=0}^{N-1}$ and can be easily implemented in optimization software. These constraints are used in all the numerical simulations presented in the following.

\begin{rem}
	In view of the previous observations when variables $(p_i)_{i=0}^{N-1}$ verify the constraints \eqref{eq:new-conv-constraint} the discrete shape constructed using \eqref{eq:formula-tri} and \eqref{eq:new-discretization} is convex. Therefore, even though the ideas regarding this discretization used the first and second derivatives of the support function, the resulting discrete framework does not depend on the regularity of the support function $p$.
	\label{rem:validity-discrete}
\end{rem}

\subsection{Convex geometry aspects.}

In the following we show further geometric arguments motivating the discretization \eqref{eq:new-discretization} and the discrete convexity constraints \eqref{eq:new-conv-constraint}. 

Let us fix the following notations:
\begin{itemize}
	\item $X(r,\theta) = (r\cos \theta,r\sin \theta)$, $O$ is the origin, $\theta_j = 2j\pi/N$, $0\leq j \leq N-1$, $N \geq 5$, $h = 2\pi/N$.
	\item $\ell(r,\theta)$ is the line going through $X(r,\theta)$ which is orthogonal to $O X(r,\theta)$.
	\item $\alpha(r,\theta)$ is the halfplane determined by $\ell(r,\theta)$, containing the origin $O$.
	\item the Hausdorff distance between two convex bodies $K_1,K_2$ is defined by $d_H(K_1,K_2) = \max\{ \sup_{x \in K_1}d(x,K_2), \sup_{y \in K_2} d(y,K1)\}$, where $d(x,K) = \inf_{y \in K} d(x,y)$.
\end{itemize}

\begin{lemma}
	1. $\ell(p_i,\theta_i)$ intersects $\alpha(p_{i-1},\theta_{i-1})\cap \alpha(p_{i+1},\theta_{i+1})$ if and only if  $p_{i+1}+p_{i-1}-2p_i\cos h \geq 0$. 
	
	2. The point $\bo A_i = p_i \bo r_i+q_i \bo t_i$ belongs to $\alpha(p_{i-1},\theta_{i-1})\cap \alpha(p_{i+1},\theta_{i+1})$ for every $p_i$ such that $p_{i+1}+p_{i-1}-2p_i\cos h \geq 0$ if and only if $q_i =(p_{i+1}-p_{i-1})/(2\sin h)$.
	\label{lem:geom}
\end{lemma}

\begin{figure}
	\centering
	\includegraphics[width=0.5\textwidth]{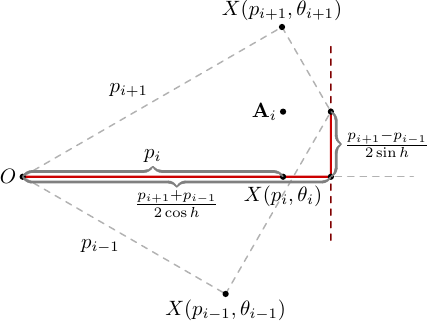}
	\caption{Geometric justification of the convexity constraint and of the discretization.}
	\label{fig:lemma-pi}
\end{figure}

\emph{Proof:} Without loss of generality, suppose that $p_{i-1},p_i,p_{i+1}>0$. The proof is self explanatory looking at the Figure \ref{fig:lemma-pi}. We sketch the main lines below.

1. Working in radial and tangential coordinates with respect to the direction $\theta_i$, the point $\ell(p_{i-1},\theta_{i-1})\cap \ell(p_{i+1},\theta_{i+1})$ has the radial coordinate $(p_{i+1}+p_{i-1})/(2\cos h)$. The statement follows.

2. The only tangential coordinate which guarantees that the point $\bo A_i$ belongs to $\alpha(p_{i-1},\theta_{i-1})\cap \alpha(p_{i+1},\theta_{i+1})$ is the tangential coordinate of the intersection $\ell(p_{i-1},\theta_{i-1})\cap \ell(p_{i+1},\theta_{i+1})$, which is exactly $(p_{i+1}-p_{i-1})/(2\sin h)$.  \hfill $\square$

\begin{rem}
	The result of Lemma \ref{lem:geom} gives a geometric motivation for the choice of discretization \eqref{eq:new-discretization} and the convexity constraint \eqref{eq:new-conv-constraint}. 
\end{rem}

As a direct consequence we have the following: 

\begin{prop}
	Let $K$ be a convex body with support function $p:[0,2\pi] \to \Bbb{R}$. For $N\geq 5$, denote by $p_j = p(2j\pi/N)$ the samples of the support function at $\theta_j$, $0\leq j \leq N-1$. Then $(p_j)_{j=0}^{N-1}$ verify inequalities \eqref{eq:new-conv-constraint}.
	\label{prop:reverse}
\end{prop}

\emph{Proof:} The result is straightforward, following Lemma \ref{lem:geom}. Notice that $\partial K \cap \ell(p_i,\theta_i)\neq \emptyset$ and $K \subset \alpha(p_{i-1},\theta_{i-1})\cap\alpha(p_{i+1},\theta_{i+1})$. \hfill $\square$

Next, we show that the set of discrete convex shapes obtained with the proposed discretization is dense in the class of convex sets. For $N \geq 5$, denote by $\mathcal K_N$ the following class of convex polygons:

\begin{multline*}
\mathcal K_N =\{\bo A_0...\bo A_{N-1} : \bo A_j = p_j \bo r_j + \frac{p_{j+1}-p_{j-1}}{2\sin h} \bo t_j,\\
 p_{j+1}+p_{j-1}-2p_j \cos \frac{2\pi}{N} \geq 0, \ \forall j=0,...,N-1\}.
\end{multline*}

\begin{thm}
	Let $K$ be an arbitrary convex body. Then for every $\varepsilon>0$ there exists $N$ large enough and a convex polygon $P_N \in \mathcal K_N$ such that $d_H(K,P_N)<\varepsilon$.
	\label{thm:approx}
\end{thm}

\emph{Proof:} Without loss of generality, suppose the origin is inside $K$. Consider $p_j = p(\theta_j)>0$, $j=0,...,N-1$. Proposition \ref{prop:reverse} implies that $(p_j)_{j=0}^{N-1}$ verify \eqref{eq:new-conv-constraint}. Denote by $Q_N$ the polygon $\bigcap_{j=0}^{N-1} \alpha(p_j,\theta_j)$. Then obviously $Q_N$ is convex and $K\subset Q_N$. 
Pick $X_j \in \ell(p_j,\theta_j)\cap K, X_{j+1} \in \ell(p_{j+1},\theta_{j+1})\cap K$ and denote $Y_j \in  \ell(p_j,\theta_j)\cap  \ell(p_{j+1},\theta_{j+1})$. The angle $X_jY_jX_{j+1}$ is equal to $\pi-h$, therefore the distance from $Y_j$ to $[X_jX_{j+1}] \subset K$ is at most $\frac{1}{2}X_jX_{j+1}\tan (h/2)$. This implies that $d_H(Q_N,K) \leq \frac{1}{2}\text{diam}(K) \tan (h/2)$. 

By construction, the polygon $\bo A_0...\bo A_{N-1}$ associated to $(p_j)_{j=0}^{N-1}$ verifies $\bo A_i \in \ell(p_i,\theta_i)$ and $\bo A_0...\bo A_{N-1} \subset Q_N$.
 More precisely, using the previous notations, we have $\bo A_j \in [Y_{j-1}Y_j]$, implying that $\bo A_0...\bo A_{N-1}$ is convex.. Therefore, 
 \[\max_{j=0,...,N-1} d(\bo A_i,K) \leq \max_{j=0,...,N-1} d(Y_j,K) \leq \frac{1}{2}\text{diam}(K)\tan (h/2).\] As a consequence we have
\[ d(\bo A_0... \bo A_{N-1},K) \leq \frac{1}{2}\text{diam}(K)\tan (h/2).\]
As $N\to \infty$ we have $h\to 0$, therefore the conclusion follows.
\hfill $\square$

In the numerical applications a large enough $N$ is chosen and discrete convex shapes are constructed following the definition of $\mathcal K_N$. In particular, the set $\mathcal K_N$ is characterized by a family of $N$ linear inequality constraints. For $N$ large enough, Theorem \ref{thm:approx} implies that any convex shape can be sufficiently well approximated using a polygon in $\mathcal K_N$.

\subsection{Width and inclusion constraints.} The second type of constraints that are of interest in this work are related to the {\bf width} or the {\bf diameter} of the convex set $K$. Such constraints can easily be formulated in terms of the support function using the quantity $w(\theta) = p(\theta)+p(\theta+\pi)$. Geometrically, $w(\theta)$ measures the distance between the two supporitng lines to $K$ orthogonal to the direction given by $\theta$. From a discrete point of view we consider $N$ even so that for any $\theta_i$ in the discretization, the antipodal point is also present $\theta_{i+N/2} = \theta_i+\pi$. In practice upper or lower bounds on width constraints can be imposed using 
\[ w_i \leq p_i+p_{i+N/2} \leq W_i,\ \ i = 0,...,N/2-1.\]

Let us list some particular cases of interest:
\begin{itemize}[noitemsep,topsep=0pt]
	\item $w_i=W_i=\overline w>0$ for all $i=0,...,N/2-1$ gives the discrete \emph{constant width constraint}. 
	\item $W_i = d$, $w_0 = d$ imposes a \emph{diameter constraint}. An upper bound on the width is considered for every direction and a lower bound is imposed for exactly two antipodal points. In this way the diameter of the set $K$ is fixed.
\end{itemize}

{\bf Inclusion} constraints can be achieved by imposing pointwise inequality constraints on the discrete values $p_i$ of the support function. Indeed, if the set $K_0$ has the support function $p_0$ the inclusion constraint $K \subset K_0$ simply reads $p \leq p_0$ on $[0,2\pi]$. In cases where $K_0$ is a polygon, it suffices to impose a finite number of constraints for orientations $\theta_i$ that are orthogonal to the edges of $K_0$.

\section{Parametric gradient of the objective function}
\label{sec:gradient}

To optimize numerically a function in an efficient way, gradient information should be used whenever available. In this section we detail the computation of the discrete gradient for various functionals used in the applications part. Given a shape $K$ and a Lipschitz vector field $V \in W^{1,\infty}(\Bbb{R}^2, \Bbb{R}^2)$ we consider the perturbed shape $(I+V)(K)$. We say that a shape functional $J$ is shape differentiable if the following expansion holds
\[ J((I+V)(K))= J(K)+J'(K)(V) + o(\|V\|_{W^{1,\infty}}).\]
For more details the classical references \cite{henrot-pierre-english}, \cite{delfour-zolesio}, \cite{sokolowski-zolesio} should be consulted. Moreover, under the assumption that the shape $K$ is convex, in all cases handled in the following, it is possible to write the shape derivative in the form 
\begin{equation}
 J'(K)(V) = \int_{\partial K} f \ V\cdot \bo n\ d\sigma ,
 \label{eq:sh-deriv-boundary}
\end{equation}
where $f$ is an integrable function on $\partial K$ and $V\cdot \bo n$ is the normal component of the perturbation vector. Alternative volume integral expressions for the shape derivatives can be given. The choice to work with boundary integrals in this work is further motivated in Remark \ref{rem:comparison-volume}.  

In the previous section the discretization of the support function using finite differences was introduced using the values $(p_i)_{i=0}^{N-1}$ at $N$ equidistant sample points in $[0,2\pi]$. In the following, we present the computation of partial derivatives of general functionals with respect to the corresponding parameters $p_i$, $i=1,...,N$. 

\subsection{General functionals.} For generic shape functionals $J(K)$, under suitable regularity assumptions which are generally valid when $K$ is convex, the corresponding shape derivative can be expressed in the form \eqref{eq:sh-deriv-boundary}. Since the discrete shape $K$ depends on the parameters $\bo (p_i)_{i=0}^{N-1}$, we can write the dependence in the form $K(\bo p)$. For $\delta \bo p \in \Bbb{R}^N$, $|\delta \bo p|\ll 1$ there is a vector field $V(\delta \bo p)$ such that $V(0)=0$ such that $K(\bo p +\delta \bo p) = (I+V(\delta \bo p))(K(\bo p))$. Assuming $V$ is differentiable at $0$ we can formally write the expansion
\begin{align*}
J(K(\bo p+\delta \bo p)) &= J((I+V(\delta \bo p))(K(\bo p))) = J((I+V(0)+DV(0) \delta \bo p)(K(\bo p))+ o(|\delta \bo p|) \\
& = J(K(\bo p))+J'(K(\bo p))(DV(0)\delta \bo p) +o(|\delta\bo p|),
\end{align*}
where $DV$ is the Jacobian matrix of $V$. 
Therefore, in order to compute the sensitivity of $J(K(\bo p))$ with respect to the parameter $p_i$ it is enough to compute the perturbation $V(\delta \bo p)$, differentiate this vector field with respect to $p_i$ and plug it in the shape derivative formula \eqref{eq:sh-deriv-boundary}.

 Note that perturbing the parameter $p_i$ with a small value $\delta p$ only changes points $\bo x_{i-1}, \bo x_i, \bo x_{i+1}$ in the discretization, in view of \eqref{eq:new-discretization}. The explicit perturbations of these points using \eqref{eq:formula-tri}, \eqref{eq:new-discretization} are given by
\[ p_i \mapsto p_i+\delta p \Longrightarrow \begin{cases}
\bo x_{i-1} &\mapsto \bo x_{i-1} +\frac{1}{\sin h} \delta p \bo t_{i-1}\\
\bo x_i &\mapsto \bo x_i+\delta p \bo r_i \\
\bo x_{i+1} &\mapsto \bo x_{i+1}-\frac{1}{\sin h}\delta p \bo t_{i+1}.
\end{cases}\]
The derivative $V_i$ of this perturbation with respect to $\delta p$ has normal components equal to $V_i \cdot \bo n = \delta_{ij}$ (the usual Kronecker delta symbol) at points $\bo x_j$, $j=0,...,N-1$. No information is known between discretization points. In view of the polygonal nature of the discrete convex shape, we make the assumption that the normal component of the perturbation vector $V_i \cdot \bo n$ is piecewise affine on the intervals $[\theta_i, \theta_{i+1}]$ corresponding to the region between boundary points $\bo x_i$ and $\bo x_{i+1}$. In order to formalize this we introduce the \emph{hat functions} $\psi_i:[0,2\pi] \to \Bbb{R}$ which are $2\pi$ periodic, continuous and piecewise affine on intervals $[\theta_i,\theta_{i+1}]$ such that $\psi_i(\theta_j) = \delta_{ij}$.

With this convention, the parametric derivative of $J$ with respect to $p_i$ becomes
\begin{equation}
\frac{\partial J(K)}{\partial p_i} = \int_{\partial K} f(\bo x) \psi_i(\theta(\bo x)) d\sigma, i=0,...,N-1,
\label{eq:grad-general}
\end{equation}
where $\theta(\bo x)$ is the orientation of the normal at the boundary point $\bo x \in \partial K$.
From a practical point of view it is necessary to transport the hat functions $\psi_i$ from $[0,2\pi]$ to $\partial K$ and perform the numerical integrations given by \eqref{eq:grad-general} for $i=0,...,N-1$. 

In the case where $\partial K$ is smooth, the parametrization \eqref{eq:param} is non-degenerate which allows us to conclude, via a change of variables, that 
\begin{equation}
\frac{\partial J(K)}{\partial p_i} = \int_0^{2\pi} f(\bo x(\theta)) \psi_i(\theta) (p(\theta)+p''(\theta)) d\theta,\ \  i=0,...,N-1.
\label{eq:grad-general-param}
\end{equation}
It can be observed that when $\bo x$ is a corner point, having multiple supporting lines, $\theta (\bo x)$ is an interval. However, a corner point corresponds to a zero curvature radius, i.e. $p(\theta)+p''(\theta)=0$.

All computations are realized using the software \texttt{FreeFEM} \cite{freefem}. Domains are meshed starting from the variables $(p_i)_{i=0}^{N-1}$, using information given by the parametrization \eqref{eq:param} with approximations \eqref{eq:new-discretization}. A discrete polygon is constructed with vertices 
\[\bo A_i = p_i \bo r_i +\frac{p_{i+1}-p_{i-1}}{2\sin h} \bo t_i, i=0,...,N-1\]
as described in Section \ref{sec:suppfunc}. If the discrete curvature radii \eqref{eq:new-conv-constraint} are non-negative, the polygon $\bo A_0...\bo A_{N-1}$ is convex. \texttt{FreeFEM} constructs the mesh starting from the polygonal line $\bo A_0...\bo A_{N-1}$. The mesh is then improved using the command \texttt{adaptmesh} in order to make it suitable for finite element computations. Finite element spaces are constructed for solving the partial differential equations involved in the computations of the objective function. In all the computations ${\bf P_2}$ finite elements are used in \texttt{FreeFEM}. The various constraints involved in the problem definitions, whose discretizations are mentioned in Section \ref{sec:suppfunc}, are formulated as linear constraints on the variables $(p_i)_{i=0}^{N-1}$ and are used in the optimization toolbox \texttt{IPOPT} \cite{IPOPT} included in \texttt{FreeFEM}. 
The sensitivity of the objective function with respect to the parameters $p_i$ is evaluated using \eqref{eq:grad-general}. The integrals are evaluated using standard \texttt{FreeFEM} routines. Figures are realized using Metapost or Matplotlib in Python.

For reproducibility purposes and in order to allow the easy adaptation of these ideas to various other problems, the codes used in the numerical simulations are available at the following repository: \href{https://github.com/bbogo/ConvexSets}{\nolinkurl{https://github.com/bbogo/ConvexSets}}

\begin{rem}
	The structure theorem for shape derivatives (see \cite[Chapter 5]{henrot-pierre-english}, \cite[Chapter 9]{delfour-zolesio}) implies that under certain regularity assumptions, shape derivatives can be written as a linear form depending on the boundary perturbation as in \eqref{eq:sh-deriv-boundary}. It is nevertheless possible to obtain shape derivatives as volume integrals. Such formulas require less regularity assumptions at the price of having derivatives on the perturbation field. While from a theoretical point of view, the two formulations (boundary vs volume integrals) are equivalent, it is no longer the case when performing numerical approximations. 
	
	In \cite{hiptmair} the authors compare the numerical errors when computing shape derivatives with the two formulations and conclude that, under additional regularity assumptions on the perturbation vector fields, the shape derivatives computed with volume integrals converge faster. A similar analysis has been performed in \cite{zhu-eigs} for the eigenvalue problems associated to the Dirichlet-Laplace eigenvalues.
	
	In this work we use boundary integrals for computing shape derivatives for multiple reasons, recalled below:
	\begin{itemize}
		\item When dealing with convex sets, the boundary integrals defining the shape derivatives are well defined for all problems under consideration. Convergence of the corresponding finite element approximations is proved in \cite{hiptmair} and \cite{zhu-eigs}.
		\item The optimization strategy presented in this work is not based on mesh perturbation techniques like in \cite{bartels_wachsmuth}. The meshed domains in the numerical computations are constructed from a set of parameters. 
		\item The perturbation fields associated to discrete perturbations in the parameters $(p_i)_{i=0}^{N-1}$ are not in $W^{2,\infty}$, so according to \cite[Remark 3.2]{hiptmair} the result \cite[Theorem 3.1]{hiptmair} does not apply.
	\end{itemize}
	\label{rem:comparison-volume}
\end{rem}

\subsection{Perimeter and Area.} In the following, the derivative of the perimeter with respect to parameters $p_i$ is investigated. Let us suppose that $K$ is a convex set, not containing segments in its boundary such that the associated support function $p$ is at least of class $C^2$. In this case, in view of the parametrization \eqref{eq:param}, the perimeter is given by integrating the arclength measure: $\Per(K) = \int_0^{2\pi} (p(\theta)+p''(\theta))d\theta$. If $p$ is of class $C^2$ we have $\int_0^{2\pi} p''(\theta)d\theta = 0$.

\begin{itemize}[noitemsep,topsep=0pt]
	\item{\bf Direct method.} The perimeter of a smooth set can be expressed in terms of the support function by the formula $\Per(K) = \int_0^{2\pi} p(\theta)d\theta$. Using the basic trapezoidal quadrature rule for intervals $[\theta_i,\theta_{i+1}]$, we have the approximation $\Per(K) \approx  \sum_{i=0}^{N-1} \frac{2\pi}{N} p_i$. For this explicit approximation formula the gradient of the perimeter with respect to $p_i$ is equal to $2\pi/N$ for every $i=0,...,N-1$. 
	\item{\bf Using the shape derivative.} For a smooth shape $K$, the shape derivative of the perimeter is given by $\Per(K)'(V) = \int_{\partial K} \mathcal H\  V\cdot \bo n$ where $\mathcal H$ is the mean curvature of $\partial K$ (equal to the curvature in dimension two). Considering, as recalled previously, the vector field $V$ corresponding to perturbation of a single variable $p_i$ in the discretization we obtain
	\begin{equation} \frac{\partial \Per(K)}{\partial p_i} = \int_0^{2\pi} \mathcal H(\bo x(\theta)) \psi_i(\theta) (p(\theta)+p''(\theta)) d\theta = \frac{2\pi}{N},
	\label{eq:grad-perim}
	\end{equation}
	where we used a change of variable and the fact that $\mathcal H(\bo x(\theta)) = 1/(p(\theta)+p''(\theta))$ at regular points $\theta \in [0,2\pi]$.
\end{itemize}

Let us now consider the non-smooth case. Suppose that segments $S_1,...,S_l$ in the boundary of $K$ exist and correspond to angles $\theta_{i_1}<\theta_{i_2}<...<\theta_{i_l} \in [0,2\pi]$. Then the same formula gives
\[ \Per(K) = \int_0^{2\pi} p(\theta)d\theta + \sum_{j=1}^{l}\left( p'(\theta_{i_j}+)-p'(\theta_{i_j}-)\right),\]
corresponding to the length of the smooth parts and the sum of the lengths of all the segments $S_j,\ j=1,...,l$. In the numerical computations, the integral of $p$ is approximated using a quadrature rule. The contribution of $p'(\theta_i+)-p'(\theta_i-)$ is added to the objective function and the gradient as soon as this difference exceeds a certain threshold, indicating a singularity.

The derivative of the area functional with respect to the parameter $p_i$ is computed by taking $f\equiv 1$ in \eqref{eq:grad-general}. Alternatively, if the support function is of class $C^1$, like in the case of shapes of constant width, the parametric derivative may be written as
\begin{equation} \frac{\partial|K| }{\partial p_i} =\int_0^{2\pi} \psi_i(\theta) (p(\theta)+p''(\theta)).
\label{eq:deriv-area-regular}
\end{equation}
In practice, for $h$ small, an approximation of \eqref{eq:deriv-area-regular} is given by $\varrho_i$, the discrete curvature radius given in \eqref{eq:new-conv-constraint}.

\section{Applications using the support function}
\label{sec:applic-supp}
In the following we present a few applications which illustrate the numerical method proposed in the previous sections. For all problems considered, the existence of solutions is discussed briefly. The proof of existence is usually straightforward, using various results recalled in the Lemma \ref{lem:existence} below. When needed, additional references containing more details are provided. The appropriate notion of convergence for convex sets is the convergence in the Hausdorff distance. Precise definitions and main properties of this set distance are found in \cite[Chapter 2]{henrot-pierre-english} or \cite[Section 1.8]{schneider}. The numerical simulations are preformed using the discretization of the support function described in Section \ref{sec:suppfunc} together with the gradient expressions described in Section \ref{sec:gradient}. In all computations below the convexity constraint is imposed using the discrete inequalities \eqref{eq:new-conv-constraint}.


We recall below some results which allow to prove existence of solutions.

\begin{lemma}	
	
1. ({Blaschke's selection theorem.}) Given a sequence $\{K_n\}$ of closed convex sets contained in a bounded set, there exists a subsequence which converges to a closed convex set $K$ in the Hausdorff metric. \cite[Theorem 1.8.7]{schneider}

2. Convexity is preserved by the Hausdorff convergence. (see \cite[p. 35]{henrot-pierre-english})

3. If $\{K_n\}$ is a sequence of non-empty closed convex sets contained in a bounded set then the Hausdorff convergence of $K_n$ to $K$ is equivalent to the uniform convergence of the support functions $p_{K_n}$ to $p$ on $\Bbb{S}^{d-1}$. (see \cite[Lemma 1.8.14]{schneider})

4. Suppose that the sequence of convex sets $\{K_n\}$ converges to the convex set $K$ in the Hausdorff topology and that $K$ has non-void interior. Then $\chi_{K_n}$ converges to $\chi_K$ in $L^1$, $|K_n| \to |K|$ and $\mathcal \Per (\partial K_n) \to \mathcal \Per (K)$ as $n \to \infty$. (see \cite[Prop 2.4.3]{bucurbuttazzo})

5. 	If $K_n$ are convex and converge to $K$ in the Hausdorff metric then $K_n$ $\gamma$-converges to $K$ and, in particular the eigenvalues of the Dirichlet-Laplace operator are continuous: $\lambda_k(K_n) \to \lambda_k(K)$. (see \cite[p. 33]{henrot-pierre-english})

6. Inclusion is stable for the Hausdorff convergence: if $\Omega$ is closed, $K_n \subset \Omega$, $K_n \to K$ implies $K\subset \Omega$. (see \cite[p. 33]{henrot-pierre-english})

7. The diameter and width constraints are continuous with respect to the Hausdorff convergence of closed convex sets. In particular if the sequence of closed convex sets $\{K_n\}$ converges to $K$ in the Hausdorff metric and each $K_n$ is of constant width $w$ then $K$ is also of constant width $w$. (a direct consequence of point 3. above)

8. The perimeter of convex sets is monotone with respect to set inclusion as shown in \cite[Lemma 2.2.2]{bucurbuttazzo}.
\label{lem:existence}
\end{lemma} 

\subsection{Area and perimeter functionals}

As a first example, the minimization of the area of a two dimensional set with minimal width $w$ is considered. In \cite[Problem 6-2]{yaglom-boltjanskii} it is proven that the solution of this problem in the class of convex sets is the equilateral triangle. This problem is used as a benchmark for the numerical algorithm since its solution in the class of convex sets is known.

\begin{problem}
	Minimize the area of a two dimensional shape under minimal width constraint: 
	\[\min \{ |\omega| : \omega  \text{ convex, having minimal width } w\}\]
	\label{prob:area-minw}
\end{problem}

\begin{figure}
	\centering 
	\includegraphics[height=0.3\textwidth]{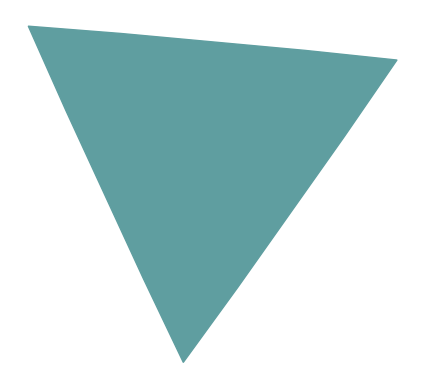}\quad
	\includegraphics[height=0.3\textwidth]{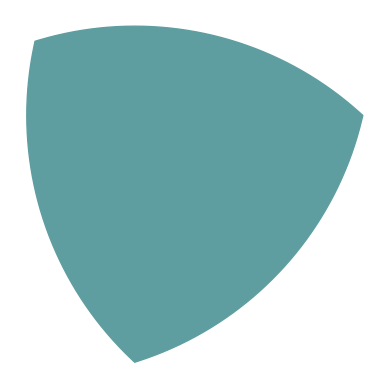}
	\caption{Minimizing the area under minimal width and constant width constraints}
	\label{fig:area-width}
\end{figure}

 The minimal width constraint is modeled numerically by considering an even number $N$ of angles  in the discretization of the support function and by imposing $p_i+p_{i+N/2} \geq w$ for $i=0,...,N/2-1$. The area and its gradient are computed using the formulas \eqref{eq:grad-general}. 
 The result given by the optimization algorithm is the equilateral triangle. The shape found by the optimization algorithm is shown in Figure \ref{fig:area-width} (left). For this problem $N=240$ discretization points were used. 

Another classical example of shape optimization problem in convex geometry for which the solution is known is the minimization of the area under constant width constraint. It is known that the solution to this problem is the Reuleaux triangle and a proof of this fact can be found in \cite[Problem 7-20]{yaglom-boltjanskii}. Following the discussion in Section \ref{sec:suppfunc}, the support function of a shape of constant width does not have discontinuities in the first derivative. Therefore, the proposed discretization is expected to handle this case without any difficulty.
\begin{problem} Minimize the area of a two dimensional convex shape of constant width $w$:
	\[\min \{ |\omega| : \omega \text{ is convex with constant width } w \}.\]
\end{problem}

 For an even number $N$ of discretization points, he discrete constant width constraint can be written in the form $p_i+p_{i+N/2} = w$ for $i=0,...,N/2-1$. For $N=240$, the result given by the numerical optimization algorithm is an approximation of the Reuleaux triangle shown in Figure \ref{fig:area-width}. Repeating the simulation using the gradient formulas for the area given by \eqref{eq:deriv-area-regular} gives similar results.

The following problem is inspired from \cite{LNconv09}. The competition between area and perimeter implies that for some values of $\mu>0$ optimal shapes are polygons. 

\begin{problem}
	Given $\mu>0$ find solutions of 
	\[\min\{ \mu |\omega|-\Per(\omega) : \omega \text{ convex }, \text{diam}(\omega)=1 \}.\]  
	\label{pb:area-perim}
\end{problem}

Existence of solutions to this problem for $\mu>0$ follows from the Blaschke selection theorem and the fact that a convex set of diameter $1$ has an upper bound on the perimeter. See Lemma \ref{lem:existence}, point 8.

\begin{figure}
	\centering
	\includegraphics[height=0.3\textwidth]{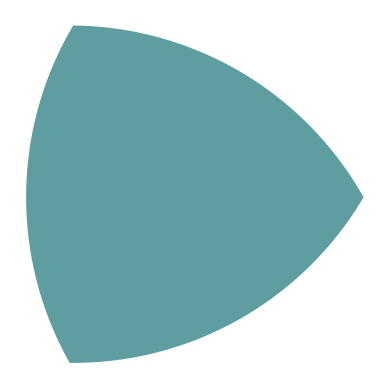}\quad
	\includegraphics[height=0.3\textwidth]{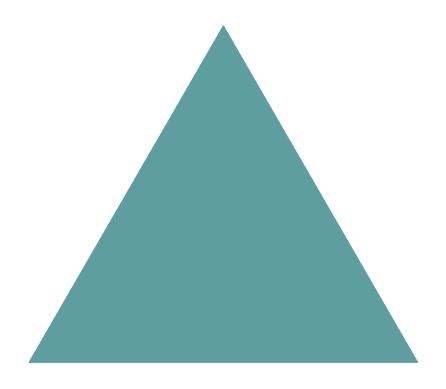}
	\caption{Numerical approximations of solutions for Problem \ref{pb:area-perim} for $\mu=0.5,1$.}
	\label{fig:area-perim}
\end{figure}

Following the choice of the parameter $\mu$, the solution changes. For $\mu>0$ small enough the solution is the Reuleaux triangle, maximizing the perimeter for a given diameter and minimizing the area. For $\mu$ large enough, the solution is a segment. For intermediary $\mu$, results of \cite{LNconv09} imply that solutions are polygons. Implementing a numerical algorithm for solving the problem is straightforward following ideas in Section \ref{sec:suppfunc}. In particular, the upper bounds on the diameter are imposed for all pairs of antipodal points, while the lower bound is imposed for one pair of antipodal points. In Figure \ref{fig:area-perim} results are shown for $\mu \in \{0.5,1\}$.

\subsection{Dirichlet Laplace eigenvalues}

Let us recall that for a Lipschitz domain $\omega$ the eigenvalues of the Laplace operator with Dirichlet boundary conditions are defined by the equation
\[ \left\{ \begin{array}{rcll}
-\Delta u & = & \lambda u & \text{ in } \omega \\
u & = & 0 & \text{ on } \partial \omega
\end{array}\right.\]
These eigenvalues form an increasing sequence $0<\lambda_1(\omega) \leq \lambda_2(\omega) \leq .. \leq \lambda_k(\omega) \to \infty$. When the shape $\omega$ is convex the first eigenvalue is simple, therefore we have $\lambda_1(\omega)<\lambda_2(\omega)$. Also it is classical that the Dirichlet-Laplace eigenvalues $\lambda_k$ are decreasing with respect to set inclusion: $\omega_1\subset \omega_2$ implies $\lambda_k(\omega_1) \geq \lambda_k(\omega_2)$. Furthermore, the behavior of the eigenvalues is well known for scalings: $\lambda_k(t\omega) = \lambda_k(\omega)/t^2$. See \cite{henrot-pierre-english} for more details. Using the scaling property one may note that minimizing $\lambda_k(\omega)$ under area constraint is equivalent to minimizing $\lambda_k(\omega)|\omega|$ without any constraints. Moreover, minimizing $\lambda_k(\omega)+|\omega|$ we obtain optimal shapes that are equivalent up to homotheties to the previous formulations.

An already classical optimization problem related to the Dirichlet Laplace eigenvalues and the convexity constraint is the minimization of the eigenvalues under area and convexity constraints. In particular, the minimization of the second eigenvalue was studied in detail in \cite{oudeteigs}, \cite{henrot_oudet} and \cite{antunes_henrot}. Therefore we formulate the following:

\begin{problem}
	Minimize $\lambda_k(\omega)|\omega|$ among convex sets.
	\label{pb:min-lamk}
\end{problem}

\begin{figure}
	\centering
	\begin{tabular}{cccc}
		\includegraphics[height=0.15\textwidth]{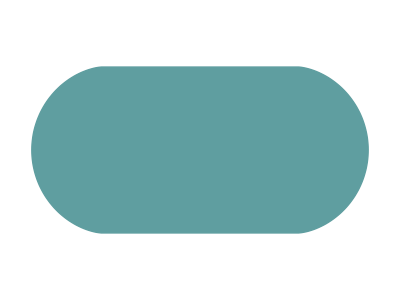}&
		\includegraphics[height=0.15\textwidth]{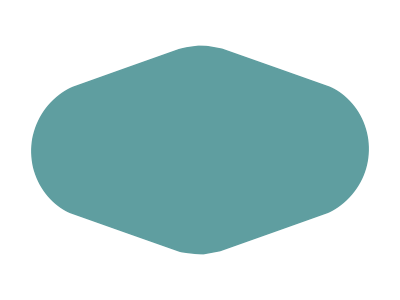}&
		\includegraphics[height=0.15\textwidth]{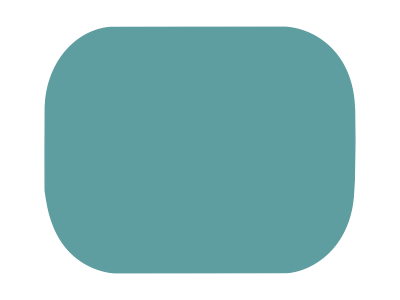}&
		\includegraphics[height=0.15\textwidth]{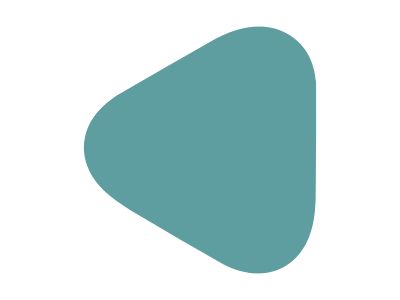}\\
		$\lambda_2=37.9855$ & $\lambda_4=65.2254$ & $\lambda_5=79.6561$ & $\lambda_6=88.5336$ \\
		\includegraphics[height=0.15\textwidth]{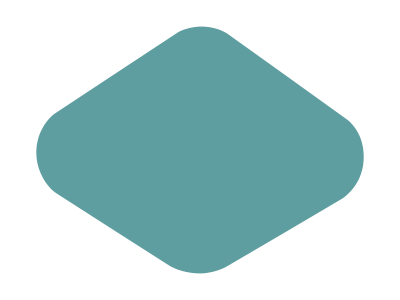}&
		\includegraphics[height=0.15\textwidth]{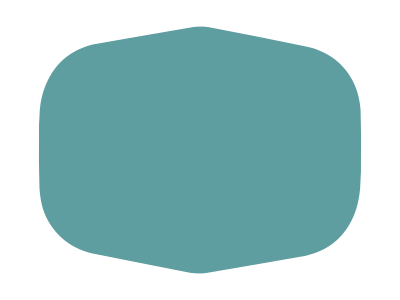}&
		\includegraphics[height=0.15\textwidth]{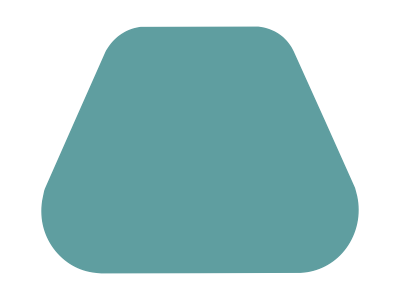}&
		\includegraphics[height=0.15\textwidth]{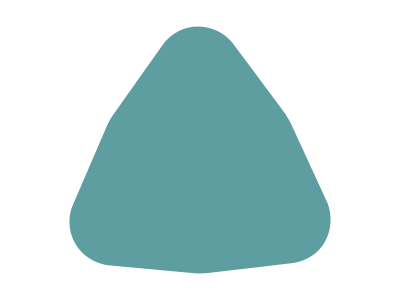}\\
	     $\lambda_7=109.1017$ & $\lambda_8=119.2929$ &
		$\lambda_9=134.9261$ & $\lambda_{10}=142.9126$ 
	\end{tabular}
	\caption{Minimization of the eigenvalues of the Dirichlet-Laplace operator under convexity and volume constraints in dimension two.}
	\label{fig:dirichlet-convex}
\end{figure}

Existence of solutions for problem \eqref{pb:min-lamk} is proved in \cite{henroteigs}. From the numerical point of view, the convexity constraint is handled using the support function as shown previously in Section \ref{sec:suppfunc}. The gradient of the eigenvalues is computed with the formula \eqref{eq:grad-general} keeping in mind that the shape derivative of a simple eigenvalue is given by 
\[ \lambda_k'(K)(V) = \int_{\partial K} |\nabla u_k|^2 V.n\  d\sigma\]
where $u_k$ is the $L^2$ normalized eigenfunction associated to $\lambda_k(\omega)$. When the eigenvalue $\lambda_k$ is multiple, the shape derivative may not exist. However, when performing numerical computations, eigenvalues are almost never multiple. Choosing the largest eigenvalue from the approximate multiplicity cluster and using it in the shape derivative formula is enough for our purposes.
The above shape derivative formula is well defined for convex sets since the corresponding eigenfunctions are in $H^2(K)$.
The results of the numerical minimization process using $N\in \{120,180\}$ are shown in Figure \ref{fig:dirichlet-convex}, together with the optimal numerical value. It can be noted that the optimal shape for $k=2$ presented here is comparable to the one obtained in \cite{antunes_henrot} and that the segments in the boundary are well captured by the parametrization proposed here. In general, the values of the objective function obtained with the current method are better than those in \cite{ABcw} since segments in the boundary are better captured. The minimization of the third eigenvalue gives the disk even without the convexity constraint as shown in  \cite{oudeteigs}, \cite{antunesf-vol}.

The minimization of the Dirichlet Laplace eigenvalues under diameter constraint was considered in \cite{BHL17}. It can be proved that when restricting ourselves to the case of constant width shapes, the maximization of these eigenvalues also makes sense. 
For more theoretical aspects regarding the existence of solutions we refer to \cite{ABcw,BHL17}. 

\begin{problem}
	Maximize $\lambda_k(K)$ under constant width constraint:
	\begin{equation*}
	 \max\{ \lambda_k(K):  p(\theta)+p(\theta+\pi) = w,\ p+p''\geq 0\}
	 \label{eq:max-lamk-const-width}
	 \end{equation*}
	\label{pb:max-lamk-const-width} 
\end{problem}

In \cite{yaglom-boltjanskii} Exercise 7-13 it is proved that the inradius of a shape of constant width is minimized by the Reuleaux triangle. Denote this value by $r>0$. Therefore, all shapes of constant width contain a disk of radius $r$, showing that the Dirichlet-Laplace eigenvalues of shapes of constant width have a finite upper bound. Coupling this with the Blaschke selection theorem we obtain existence of solutions for Problem \ref{pb:max-lamk-const-width} for any $k\geq 1$. 

The computation of the Dirichlet-Laplace eigenvalues is realized using finite elements and the convexity and constant width constraints are imposed as indicated in Section \ref{sec:suppfunc}. Numerical simulations indicate that the Reuleaux triangle is the solution to Problem \ref{pb:max-lamk-const-width} for $1\leq k \leq 10$.


Sets of given minimal width $w$ have an lower bound on the inradius \cite[Exercise 6-2]{yaglom-boltjanskii}. Therefore the corresponding Dirichlet-Laplace eigenvalues have an upper bound. Restricting the sets to a closed bounded ball increases the Dirichlet-Laplace eigenvalues. These considerations, together with the Blaschke selection lemma imply that the following problem has solutions.

\begin{problem}
	Maximize $\lambda_k(K)$ under minimal width constraint $w$:
	\[ \max\{ \lambda_k(K) : p(\theta)+p''(\theta+\pi) \geq w, p+p''\geq 0\}.\]
	\label{pb:max-lamk-width}
\end{problem}

Numerical simulations, using the previously described numerical framework, show that solutions to Problem \ref{pb:max-lamk-width} are equilateral triangles for $1\leq k \leq 10$.

Numerical results regarding Problems \ref{pb:max-lamk-const-width} and \ref{pb:max-lamk-width} are in accord with extremality results concerning the area functional, regarding the Reuleaux triangle in the class of shapes of constant width and the equilateral triangle in the class of shapes with minimal width. Exploiting the monotonicity of the eigenvalues with respect to inclusions might lead to a theoretical proof of these new numerical conjectures. 

\subsection{General functionals}

In this section we illustrate how the numerical framework applies to the problems proposed in \cite{bartels_wachsmuth}. In the following $Q$ is a compact convex subset of $\Bbb{R}^d$. The problems of interest are PDE constrained optimization problems of the form:
\begin{equation}
  \min \left\{J (K)=  \int_K j(x,u,\nabla u) :  u \in H_0^1(K), \ -\Delta u = f \text{ in }K,\ K \text{ convex}, K\subset Q \right\},
  \label{eq:problem_bartels_wachsmuth}
\end{equation}
where $j: Q\times \Bbb{R} \times \Bbb{R}^d \to \Bbb{R}$ satisfies suitable growth conditions. For $d=2$, in \cite[Prop. 3.1]{bartels_wachsmuth} it is proved that if $|j(x,u,v)|\leq a(x)+c(|u|^p+|v|^2)$ for $c \geq 0$ and $a \in L^1(Q)$, $p<\infty$ then problem \eqref{eq:problem_bartels_wachsmuth} has solutions. 

The theory regarding the shape derivatives of the functional appearing in \eqref{eq:problem_bartels_wachsmuth} is classical and recalled in \cite{bartels_wachsmuth}. In particular, it follows that 
\begin{equation}
J'(K)(V) = \int_{\partial \Omega} \left(j(x,u,\nabla u)-\frac{\partial u}{\partial n} \frac{\partial p}{\partial n} \right) V\cdot n \ d\sigma,
\label{eq:sh-deriv-bartels}
\end{equation}
where $u, p \in H_0^1(K)$ solve the problems $-\Delta u = f$ and the adjoint problem $-\Delta p = -j'_u(x,u,\nabla u)+\di j'_v(x,u,\nabla u)$. It is well known that for $K$ convex the solutions $u$ of the state problem and $p$ of the adjoint problem are in $H^2(K)$ \cite{grisvard}. This implies that the integral in \eqref{eq:sh-deriv-bartels} is well defined.

In \cite{bartels_wachsmuth} the particular case $j(x,u,v)=u$ was considered for the functions
\begin{equation}
f_1(x_1,x_2) = 20(x_1+0.4-x_2^2)^2+x_1^2+x_2^2-1
\label{eq:fun-BW1}
\end{equation}
and
\begin{align}
f_2(x_1,x_2) = -\frac{1}{2}+\frac{4}{5}(x_1^2+x_2)^2+2& \sum_{i=0}^{n-1} \exp(-8((x_1-y_{1,i})^2+(x_2-y_{2,i})^2))\notag \\
-&\sum_{i=0}^{n-1} \exp(-8((x_1-z_{1,i})^2+(x_2-z_{2,i})^2))
\label{eq:fun-BW2}
\end{align}
for $n=5$, $y_{1,i} = \sin((i+1/2)2\pi/n)$, $y_{2,i} = \cos((i+1/2)2\pi/n)$, $z_{1,i}= \frac{6}{5}\sin(i2\pi/n)$, $z_{2,i}=\frac{6}{5} \cos(i2\pi/n)$. This gives rise to the following:
\begin{problem}
	Solve problem \eqref{eq:problem_bartels_wachsmuth} for $J(K) = \int_K u \ dx$ and $f=f_i$, $i=1,2$ given in \eqref{eq:fun-BW1} and \eqref{eq:fun-BW2}.
	\label{pb:BW}
\end{problem}

\begin{figure}
	\centering
	\begin{tabular}{cc}
	\includegraphics[height=0.4\textwidth]{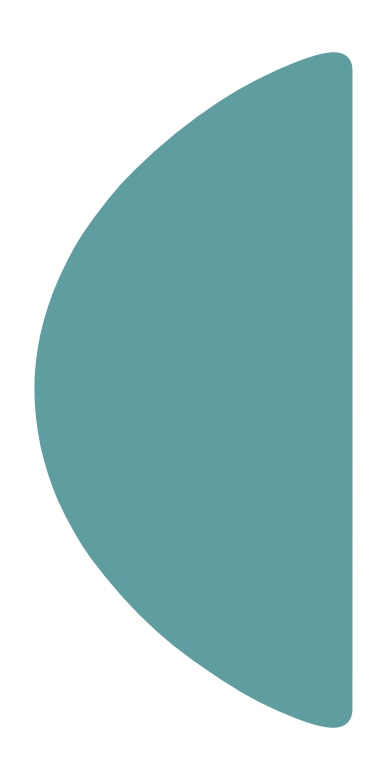}&
	\includegraphics[height=0.4\textwidth]{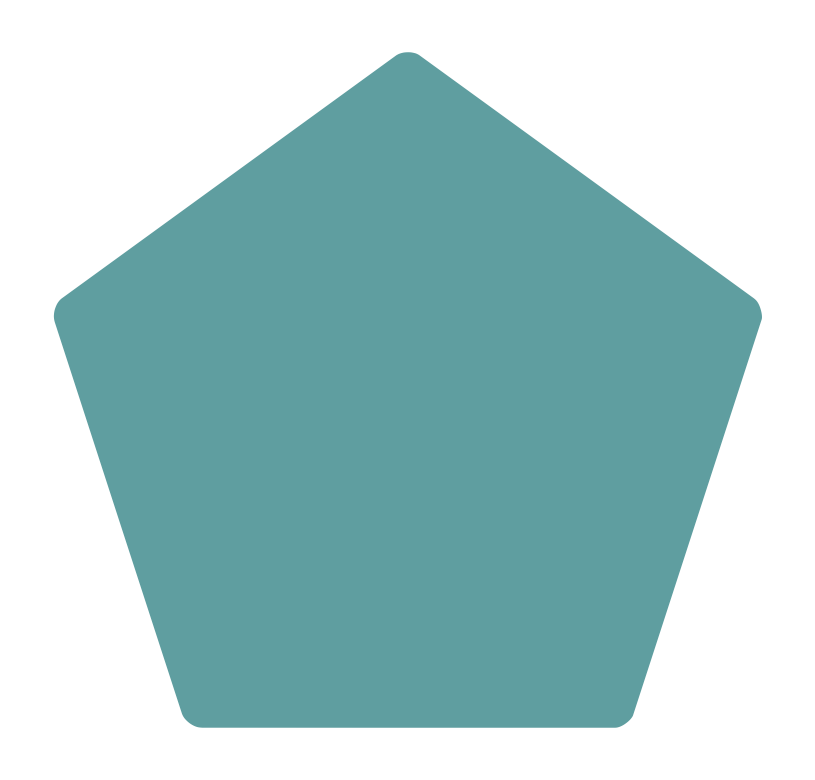}\\
	\includegraphics[height=0.4\textwidth]{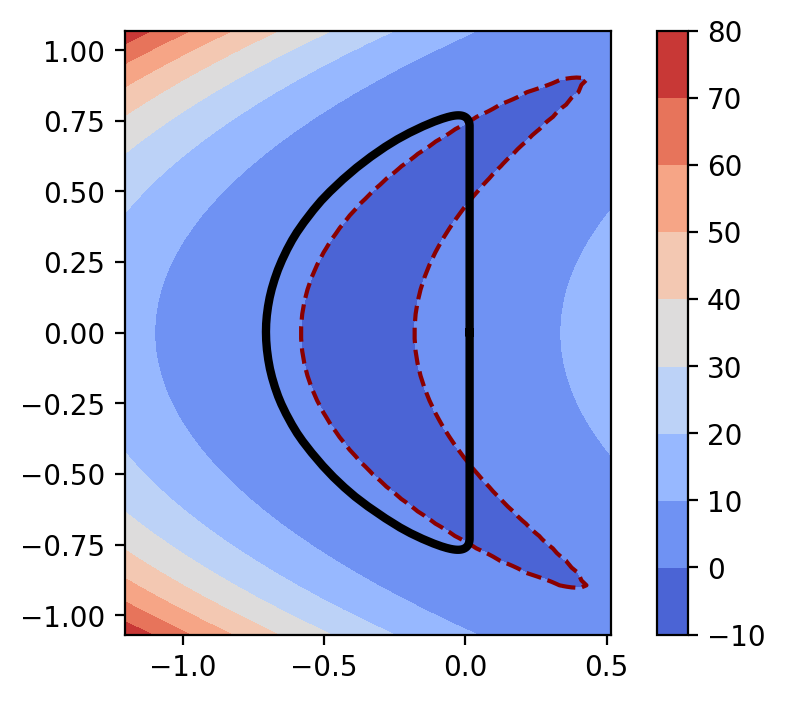}&
	\includegraphics[height=0.4\textwidth]{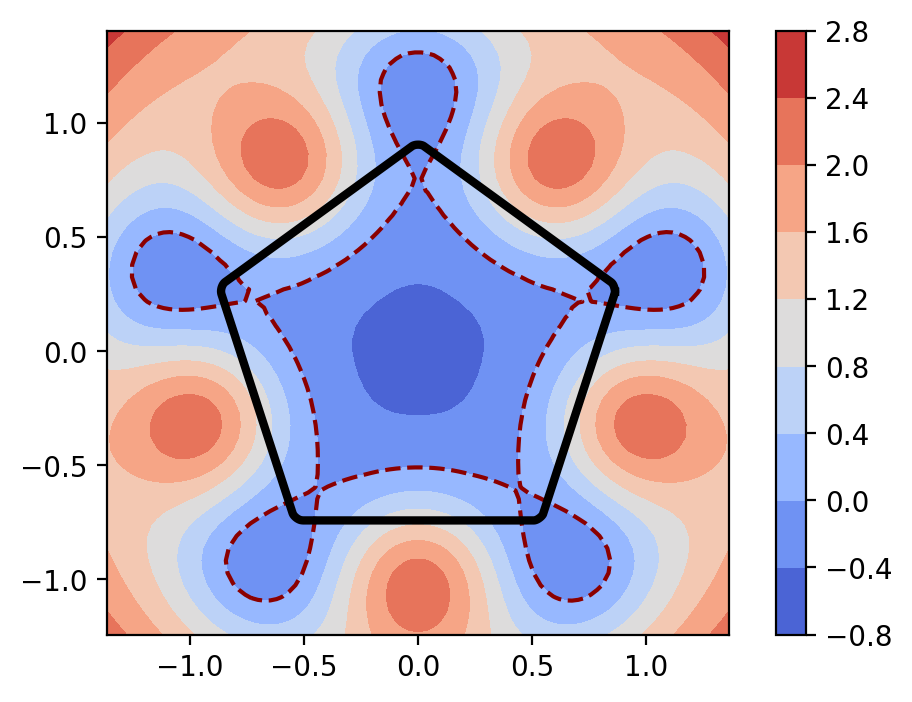}\\
	$J(K^*) = $ \texttt{-6.12084e-3} & $J(K^*) =$ \texttt{-2.76243e-2}
	\end{tabular}
	\caption{Solutions for Problem \ref{pb:BW}. The optimal shapes $K^*$ are represented and are superposed with the corresponding function $f_i$, $i=1,2$ defined in \eqref{eq:fun-BW1}, \eqref{eq:fun-BW2}.}
	\label{fig:BW}
\end{figure}

 Using the proposed numerical framework, it is straightforward to solve this problem numerically. The finite element problems are solved using ${\bf P_2}$ finite elements in \texttt{FreeFEM}. The resulting shapes together with the associated numerical optimal values are shown in Figure \ref{fig:BW}. It can be observed that the functions $f_1, f_2$ are constructed such that the sets $\{f_i\leq 0\}$ are non convex. Minimizing the objective function forces $K$ to be close to the sets $\{f_i\leq 0\}$. On the other hand, the convexity constraint imposed on $K$ is an obstacle for this, which forces parts of the optimal sets $K^*$ to be segments.

\section{Alternative discretization: the Gauge function}

\label{sec:gauge}

A convex shape $K$ with non-void interior is well characterized using a radial function $\rho : [0,2\pi] \to (0,+\infty)$ with respect to an interior point $O$. The radial function verifies $\rho(\theta) = |OX_\theta|$ where $X_\theta \in \partial K$ is the intersection of the line through $O$ having direction $(\cos \theta,\sin \theta)$ with $\partial K$. Given a radial function $\rho(\theta)$ which is of class $C^2$ at least, the curvature of $K$ for the radial coordinate $\theta$ is given by
\[ \kappa(\theta) = \frac{\rho^2(\theta)+2(\rho'(\theta))^2-\rho(\theta)\rho''(\theta)}{(\rho(\theta)^2+(\rho'(\theta))^2)^{3/2}}.\]
It can be readily checked that using the gauge function, defined by $\gamma : [0,2\pi] \to (0,+\infty)$, $\gamma(\theta)=1/\rho(\theta)$ the sign of the curvature $\kappa (\theta)$ is given by the sign of $\gamma+\gamma''$. In other words, if $\gamma$ is of class $C^2$ then $\gamma$ is the gauge function of a convex set if and only if 
\begin{equation} \gamma(\theta)+\gamma''(\theta) \geq 0, \text{ for every } \theta \in [0,2\pi].
\label{eq:convexity-gauge}
\end{equation}
As recalled in the introduction, the gauge function of a convex body is the support function of the polar body $K^\circ = \{y \in \Bbb{R}^d : x\cdot y \leq 1, \forall x \in K\}$.

As in the case of the support function, described in Section \ref{sec:suppfunc} we consider a discretization $\theta_j = jh$, $0\leq j \leq N-1$, with $h = 2\pi/N$. The values of the gauge function at the points $\theta_i$ are approximated by $\gamma_i \approx \gamma(\theta_i)$. Note that by definition we have $\gamma_i>0$. The discretization of the convexity constraint \eqref{eq:convexity-gauge} using centered finite differences gives
\begin{equation}
 \gamma_i+\frac{\gamma_{i+1}-2\gamma_i+\gamma_{i-1}}{h^2} \geq 0, \text{ for every } 0 \leq i \leq N-1.
 \label{eq:discrete-gauge-fd}
\end{equation}

On the other hand, for three consecutive angles $\theta_{i-1},\theta_i,\theta_{i+1}$ we may consider the triangle with vertices $\bo A_i = (1/\gamma_i) \bo r_i$, with $\bo r_i$ given in \eqref{eq:radial-tangent} and computing its oriented area using \eqref{eq:area-tri} gives 
\begin{equation}
\text{Area}(\Delta \bo A_{i-1} \bo A_i\bo A_{i+1}) = \frac{(\gamma_{i-1}+\gamma_{i+1}-2\gamma_i\cos h)\sin h }{2\gamma_{i-1}\gamma_i\gamma_{i+1}}.
\label{eq:area-gauge}
\end{equation} 
The Mathematica script performing the symbolic computation is given in the Appendix.
This implies that the rigorous convexity condition from the discrete point of view is
\begin{equation}
\gamma_{i-1}+\gamma_{i+1}-2\gamma_i\cos h \geq 0, \text{ for every } 0 \leq i \leq N-1.
\label{eq:convexity-gauge-rigorous}
\end{equation}
In view of the equality $2\cos h = 2-h^2+O(h^4)$, inequalities \eqref{eq:discrete-gauge-fd} and \eqref{eq:convexity-gauge-rigorous} are equivalent up to a term of order $O(h^4)$. However, for small $h$ \eqref{eq:discrete-gauge-fd} is a consequence of \eqref{eq:convexity-gauge-rigorous}, but not the other way around. 

It can be observed that the rigorous discrete convexity constraint \eqref{eq:convexity-gauge-rigorous} is the same as the rigorous discrete convexity constraint for the support function \eqref{eq:new-conv-constraint}. Therefore, given a set of parameters $(p_i)_{i=0}^{N-1}=(\gamma_i)_{i=0}^{N-1}$, verifying the constraints \eqref{eq:convexity-gauge-rigorous}, the discrete shapes constructed using the proposed discretization for the support function and the gauge functions are both convex. The famous Mahler inequalities \cite{mahler} study the bodies that minimize or maximize the product of the volume of the body $K$ and the volume of the polar body $K^\circ$. In view of the previous observations, numerical tools can be constructed based on the support and gauge functions, which can parametrize simultaneously a convex shape $K$ and its polar $K^\circ$ using a single set of parameters. 

The aspects shown previously show that it is straightforward to implement the numerical optimization of shapes with respect to the parameters $\gamma_i$, by imposing the linear inequalities \eqref{eq:convexity-gauge-rigorous} in a numerical optimization software. A straightforward computation shows that given the shape derivative formula \eqref{eq:sh-deriv-boundary} the sensitivity of the functional $J$ with respect to the parameter $\gamma_i$ is given by
\begin{equation}
\frac{\partial J(K)}{\partial \gamma_i} = -\frac{1}{\gamma_i^2}\int_{\partial K} f(\bo x) \ \psi_i(\theta(\bo x)) (\bo n \cdot \bo r_i)d\sigma . 
\label{eq:grad-gauge}
\end{equation}
The functions $\psi_i$ are the same hat functions as the ones used in \eqref{eq:grad-general}. In \eqref{eq:grad-gauge} $\theta(\bo x)$ denotes the angle of the point $\bo x$ in radial coordinates and $\bo n$ is the corresponding normal vector.

\begin{rem}
	The characterization of discrete convex shapes using the gauge function is more straightforward compared to the support function. However, diameter or witdh constraints cannot be handled in a direct way as it was the case for support functions.
	
	Recall that support function has singularities when segments are present in the boundary. The gauge function is not singular on segments, but has singularities at corners.
	\label{rem:comparison-support-gauge}
\end{rem}

We conclude this section with a few numerical examples. As in the case of the support function, \texttt{FreeFEM} is used for solving the PDEs involved in the computation of the objective function and \texttt{IPOPT} is used for handling the optimization process and the constraints. 

\begin{figure}
	\centering
	\begin{tabular}{cc}
	\includegraphics[height=0.3\textwidth]{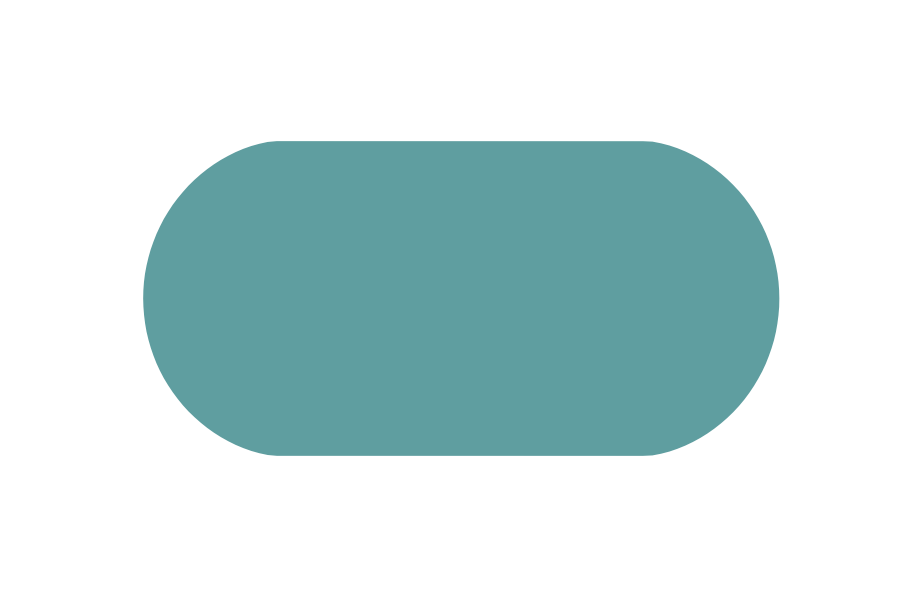}&
	\includegraphics[height=0.3\textwidth]{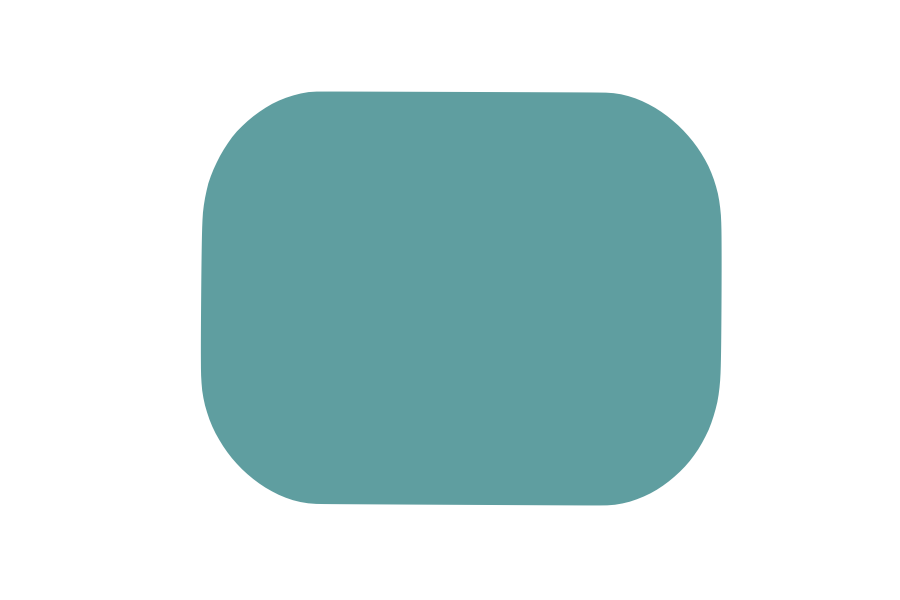}\\
	$\lambda_2 = 37.986$ & 
	$\lambda_5 = 79.644$	
\end{tabular}
	\caption{Solutions for Problem \ref{pb:min-lamk} obtained using the gauge function together with optimal numerical values at unit area.} 
	\label{fig:gauge-eigs}
\end{figure}

In Figure \ref{fig:gauge-eigs} solutions to Problem \eqref{pb:min-lamk} for $k \in \{2,5\}$ obtained using the parametrization based on the gauge function are presented. It can be observed that results are comparable with those in Figure \ref{fig:dirichlet-convex} and segments in the boundary are captured efficiently.

The recent article \cite{jimmy_2021} shows that in the class of convex sets, the maximization of the first Dirichlet-Laplace eigenvalue with inclusion constraints is well posed and any maximizing set is polygonal in the free region. However, the optimal shapes are not known in general. This motivates the numerical study the following problem.

\begin{figure}
	\centering
	\includegraphics[width=0.3\textwidth]{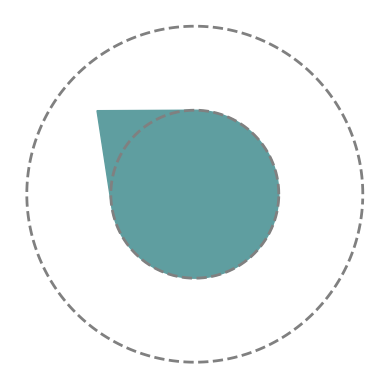}
	\includegraphics[width=0.3\textwidth]{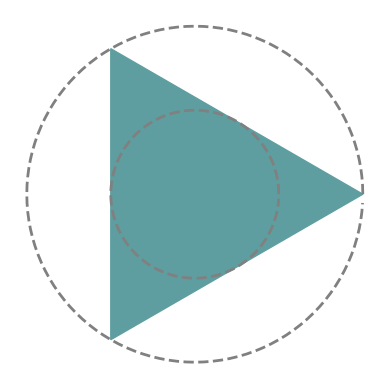}
	\caption{
	Results for Problem \ref{pb:max-lamk-inclusion} when $c \in \{\pi+0.3,3\sqrt{3}\}$ and $D_1,D_2$ are concentric disks of radii $1$ and $2$.}
	\label{fig:max-lam1-inclusions}
\end{figure}

\begin{problem}
	Given open convex sets $D_1\subset D_2$ and $c \in (|D_1|,|D_2|)$ solve
	\[ \max\{ \lambda_1(\omega) : \omega \text{ convex }, D_1 \subset \omega \subset D_2, |\omega| = c\}.\]
	\label{pb:max-lamk-inclusion}
\end{problem}
The numerical setting is strictly similar as the one used in Problem \ref{pb:min-lamk}. The discretization of the shape is realized using the gauge function. We consider the case where $D_1$ is the unit disk centered at the origin and $D_2$ is the disk of radius $2$ centered at the origin. We illustrate in Figure \ref{fig:max-lam1-inclusions} results obtained for $c\in \{ \pi+0.3, 3\sqrt{3}\}$. In particular, $c=3\sqrt{3}$ corresponds to the equilateral triangle inscribed in $D_2$ whose incircle is $D_1$. Of course, one can study in detail the behavior of the solutions with respect to the volume constraint $c$, but this is not he main purpose of this article. 

In the following we consider problems involving functionals depending on the convex body $K$ and its polar $K^\circ$, that combine the usage of the support and the gauge function. 

\begin{problem}
	Minimize $|K|\cdot |K^\circ|$ when  
	\begin{itemize}
		\item $K$ is convex and symmetric with respect to the origin.
		\item $K$ is a general convex body containing the origin.
	\end{itemize} 
\label{pb:mahler}
\end{problem}

\begin{figure}
	\centering
	\begin{tabular}{cc}
	\includegraphics[height=0.3\textwidth]{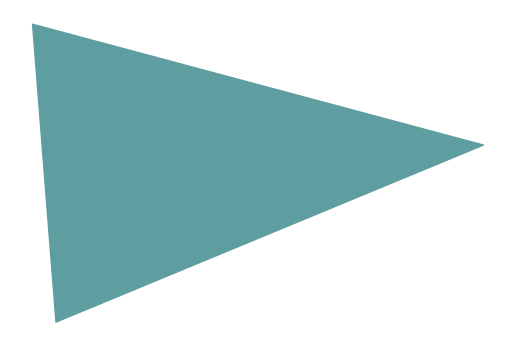}&
	\includegraphics[height=0.3\textwidth]{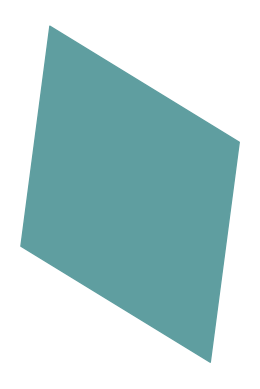}\\
	$|K|\cdot |K^\circ| = 6.750$ & $|K|\cdot |K^\circ| = 8.000$
\end{tabular}
	\caption{Numerical minimizers obtained for Problem \ref{pb:mahler}.}
	\label{fig:mahler}
\end{figure}

In view of the results shown in \cite{mahler} the solutions to the problem above are parallelograms and triangles having the centroid at the origin, respectively. In the numerical algorithm the body $K$ is parametrized using the support function as described in Section \ref{sec:suppfunc} while its polar body $K^\circ$ is characterized using the gauge function with the same parameters. The functional being scale invariant, pointwise upper and lower bounds are imposed for every variable in the parametrization to improve the stability of the optimization algorithm. The symmetry with respect to the origin is implemented by choosing an even number $N$ of equidistant parametrization angles and parameters which verify $p_i = p_{i+N/2}$, $i=0,...,N/2-1$. The results given by the numerical algorithm are given in Figure \ref{fig:mahler} and they coincide with the analytical ones discussed in \cite{mahler}. In particular, the minimization in the class of convex sets symmetric with respect to the origin gives a parallelogram, while the minimization in the class of general convex bodies containing the origin gives a triangle. 

\section{Conclusions}

This paper illustrates how the support function and the gauge function can be used to approximate solutions to shape optimization problems among convex sets. Compared to the previous works \cite{ABcw}, \cite{BHcw12}, the methods proposed here can capture well the presence of segments in the boundary. Compared to \cite{steklov_diam} a rigorous discrete convexity condition is found and implemented. The only constraint for functionals to be optimized with the proposed method is the existence of the shape derivative in the form \eqref{eq:sh-deriv-boundary}. Therefore, functionals involving solutions of partial differential equations can be efficiently handled. 

From a practical point of view the parametrizations involving the support function and the gauge function have similar complexity, notably the discrete convexity condition being the same. There are, however, some differences which we underline below:
\begin{itemize}
	\item The support function allows to easily formulate width, constant-width and diameter constraints. Inclusion constraints can be easily formulated using both parametrizations.
	\item Support functions have singularities for segments in the boundary, while gauge functions have singularities at corners (or angular points).
\end{itemize}
 While discontinuities in the derivative can be handled by the proposed method, the discretization method should be chosen to be best adapted to the problem studied. For example: if the solution is not expected to have segments in the boundary, the support function can be used; if the solution is not expected to have angular points in the boundary, the gauge function can be used. 

A wide range of applications is given, illustrating the versatility of the method for the study of two dimensional problems. Codes used for some of the problems illustrated in the article are available at \href{https://github.com/bbogo/ConvexSets}{\nolinkurl{https://github.com/bbogo/ConvexSets}}.

{\bf Acknowledgments:} The author thanks the authors of \cite{bartels_wachsmuth} for sharing information about the numerical optimizers from their work. The author was partially supported by the ANR Shapo (ANR-18-CE40-
0013) programme.

\bibliography{./DiscreteConvex}
\bibliographystyle{abbrv}

\section*{Statements and Declarations}

{\bf Funding.} The author was partially supported by the ANR Shapo (ANR-18-CE40-0013) programme.

{\bf Competing Interests.} The author has no relevant financial or non-financial interests to disclose.

\appendix

\section{Code for symbolic computations}

In order to avoid writing the tedious computations leading to the formulas \eqref{eq:expansion-area-fd}, \eqref{eq:area-tri-rigorous}, \eqref{eq:area-gauge}, scripts performing the equivalent symbolic computations in Mathematica are provided below.

Mathematica script for the computation \eqref{eq:expansion-area-fd}: 
\begin{verbatim}
p0 := (\[Rho]1 - p1)*h^2 + 2 p1 - p2
p4 := (\[Rho]3 - p3)*h^2 + 2 p3 - p2
p2 := (p1 + p3 - \[Rho]2*h^2)/(2 - h^2)
x1 := p1*Cos[t - h] - q1*Sin[t - h]
y1 := p1*Sin[t - h] + q1*Cos[t - h]
x2 := p2*Cos[t] - q2*Sin[t]
y2 := p2*Sin[t] + q2*Cos[t]
x3 := p3*Cos[t + h] - q3*Sin[t + h]
y3 := p3*Sin[t + h] + q3*Cos[t + h]
q1 := (p2 - p0)/(2*h)
q2 := (p3 - p1)/(2*h)
q3 := (p4 - p2)/(2*h)
S:= 1/2*((x2 - x1)*(y3 - y2) - (x3 - x2)*(y2 - y1))
A = Series[TrigReduce[S], {h, 0, 3}]
\end{verbatim}

Mathematica script for the computation \eqref{eq:area-tri-rigorous}: 
\begin{verbatim}
p0 := \[Rho]1*(2 - 2 Cos[h]) + 2*Cos[h]*p1 - p2
p4 := \[Rho]3*(2 - 2 Cos[h]) + 2*Cos[h]*p3 - p2
p2 := (p1 + p3 - \[Rho]2*(2 - 2 Cos[h]))/(2*Cos[h])
x1 := p1*Cos[t - h] - q1*Sin[t - h]
y1 := p1*Sin[t - h] + q1*Cos[t - h]
x2 := p2*Cos[t] - q2*Sin[t]
y2 := p2*Sin[t] + q2*Cos[t]
x3 := p3*Cos[t + h] - q3*Sin[t + h]
y3 := p3*Sin[t + h] + q3*Cos[t + h]
q1 := (p2 - p0)/(2*Sin[h])
q2 := (p3 - p1)/(2*Sin[h])
q3 := (p4 - p2)/(2*Sin[h])
S:= 1/2*( (x2 - x1)*(y3 - y2) - (x3 - x2)*(y2 - y1))
A = Simplify[TrigReduce[S]]
\end{verbatim}

Mathematica script for the computation \eqref{eq:area-gauge}: 
\begin{verbatim}
x1 := 1/\[Gamma]1*Cos[t - h]
y1 := 1/\[Gamma]1*Sin[t - h]
x2 := 1/\[Gamma]2*Cos[t]
y2 := 1/\[Gamma]2*Sin[t]
x3 := 1/\[Gamma]3*Cos[t + h]
y3 := 1/\[Gamma]3*Sin[t + h]
S:= 1/2*( (x2 - x1)*(y3 - y2) - (x3 - x2)*(y2 - y1))
A = Simplify[TrigExpand[TrigReduce[S]]]
\end{verbatim}
\end{document}